\documentclass[11pt]{article}
\usepackage{amsmath}
\usepackage{centernot}
\usepackage{listings}
\usepackage{amssymb}
\usepackage{array}
\usepackage{amsthm}
\usepackage{mathtools}
\usepackage[mathscr]{euscript}
\usepackage{relsize}
\usepackage{mathrsfs}
\usepackage{tikz}
\usetikzlibrary{matrix,arrows}
\usepackage{blindtext}
\usepackage[utf8]{inputenc}
\usepackage[margin=1.0in]{geometry}
\usepackage{enumitem}
\usepackage{url}

\theoremstyle{definition}
\newtheorem{defn}{Definition}[section]
\theoremstyle{plain}
\newtheorem{thm}[defn]{Theorem}
\newtheorem{cor}[defn]{Corollary}

\newtheorem{lm}[defn]{Lemma}
\newtheorem{pr}[defn]{Proposition}

\newtheorem{conj}[defn]{Conjecture}
\theoremstyle{remark}
\newtheorem{rem}[defn]{Remark}

\DeclareMathOperator{\Div}{Div}
\DeclareMathOperator{\Symb}{Symb}
\DeclareMathOperator{\Sym}{Sym}

\DeclareMathOperator{\Ind}{Ind}
\DeclareMathOperator{\End}{End}
\DeclareMathOperator{\Sp}{Sp}

\DeclareMathOperator{\Hom}{Hom}

\DeclareMathOperator{\ord}{ord}
\DeclareMathOperator{\Norm}{N}

\newcommand{\ds}{\displaystyle}

\newcommand{\R}{\mathbb{R}}
\newcommand{\C}{\mathbb{C}}
\newcommand{\Prj}{\mathbb{P}}

\newcommand{\D}{\mathbb{D}}

\newcommand{\Z}{\mathbb{Z}}
\newcommand{\Q}{\mathbb{Q}}
\newcommand{\T}{\mathbb{T}}
\newcommand{\We}{\mathcal{W}}

\newcommand{\p}{\mathfrak{p}}
\newcommand{\pa}{\mathfrak{P}}
\newcommand{\q}{\mathfrak{q}}
\newcommand{\ga}{\mathfrak{a}}

\newcommand{\OX}{\mathscr{O}}

\newcommand{\Gal}{\text{Gal}}

\newcommand{\GL}{\text{GL}}

\newcommand{\A}{\mathbb{A}}

\newcommand{\diagram}{\begin{tikzpicture}[description/.style={fill=white,inner sep=2pt}]
\matrix (m) [matrix of math nodes, row sep=3em,
column sep=2.5em, text height=1.5ex, text depth=0.25ex]}

\newcommand{\niceArray}{\arraycolsep=1.4pt\def\arraystretch{2.2}}
\newcommand{\fG}{\widehat{\mathbb{G}}}
\newcommand{\niceArrayPolys}{\arraycolsep=1.4pt\def\arraystretch{1.4}}

\renewenvironment{thebibliography}[1]{%
   \begin{oldthebibliography}{#1}%
     \setlength{\itemsep}{-.3ex}%
}%
{%
   \end{oldthebibliography}%
}

\newcommand{\Addresses}{{
  \bigskip
  \footnotesize

  Joseph Ferrara, \textsc{Department of Mathematics, University of California, San Diego, 9500 Gilman Drive \# 0112, La Jolla, CA 92093}\par\nopagebreak
  \textit{E-mail address:} \texttt{jferrara@ucsd.edu}

}}

\begin{document}

\setlength\abovedisplayskip{3pt plus 1pt minus 1pt}
\setlength\belowdisplayskip{3pt plus 1pt minus 1pt}

\title{A $p$-adic Stark conjecture in the rank one setting}
\author{Joseph W. Ferrara}
\date{}

\maketitle
\begin{abstract} We give a new definition of a $p$-adic $L$-function for a mixed signature character of a real quadratic field and for a nontrivial ray class character of an imaginary quadratic field. We then state a $p$-adic Stark conjecture for this $p$-adic $L$-function. We prove our conjecture in the case when $p$ is split in the imaginary quadratic field by relating our construction to Katz's $p$-adic $L$-function. We also provide numerical evidence for our conjecture in three examples.
\end{abstract}

\tableofcontents

\section{Introduction}

Let $F$ be a number field and let
$$\chi:G_F\longrightarrow\C^\times$$
be a continuous one dimensional representation of the absolute Galois group of $F$. Let $K$ be the fixed field of the kernel of $\chi$. For the rest of this article, let $p$ be an odd prime number, let $\overline{\Q}$ be an algebraic closure of $\Q$, and fix embeddings $\iota_\infty:\overline{\Q}\hookrightarrow\C$ and $\iota_p:\overline{\Q}\hookrightarrow\C_p$.

Via the Artin map, to $\chi$ we may associate the complex Hecke $L$-function, $L(\chi,s)$, defined by the series
$$L(\chi,s) = \sum_{\ga\subset\OX_F}\frac{\chi(\ga)}{\Norm\ga^s}$$
for Re$(s)>1$. The function $L(\chi,s)$ has a meromorphic continuation to the whole complex plane. In the late 1970s, in a series of papers, Stark made precise conjectures concerning the leading term of the Taylor series expansion at $s = 0$ of $L(\chi,s)$ (\cite{St2}, \cite{St3}, \cite{St4}). Stark's conjectures relate the leading term of $L(\chi,s)$ at $s=0$ to the determinant of a matrix of linear combinations of logarithms of units in $K$. His conjectures refine Dirichlet's class number formula. Stark proved his conjectures when the field $F$ is equal to $\Q$ or to an imaginary quadratic field. In general the conjectures are open.  

Around the same time that Stark made his conjectures, $p$-adic $L$-functions were constructed interpolating the critical values of complex Hecke $L$-functions for general number fields. This vastly generalized Kubota and Leopoldt's work on the $p$-adic Riemann zeta function. When $F$ is a totally real field and $\chi:G_F\rightarrow\C^\times$ is a totally even character, Cassou-Nogues (\cite{CN}), and then Deligne and Ribet (\cite{DR}) defined a $p$-adic meromorphic function
$$L_p(\chi,s):\Z_p\longrightarrow\C_p$$
determined by the following interpolation property: for all $n\in\Z_{\leq 0}$,
\begin{equation}\label{DRInt}L_p(\chi,n) = \prod_{\p\mid p}(1 - \chi\omega^{n-1}(p)\Norm\p^{-n})L(\chi\omega^{n-1}, n)\end{equation}
where $\omega$ is the Teichm\"{u}ller character. Siegel and Klingen (\cite{Sie}) showed that the values $L(\chi\omega^{n-1},n)$ lie in the field obtained by adjoining the values of $\chi\omega^{n-1}$ to $\Q$. The equality (\ref{DRInt}) takes place in  $\overline{\Q}$.

Now let $F$ be a CM field with maximal totally real subfield $E$. A prime $p$ is called ordinary for $F$ if every prime above $p$ in $E$ splits in $F$. For such primes $p$, Katz (\cite{K1},\cite{K2}) defined a $p$-adic $L$-function associated to any ray class character $\chi:G_F\rightarrow\C^\times$. Katz's $p$-adic $L$-function interpolates the values of complex $L$-functions of algebraic Hecke characters with nonzero infinity type. To specify the interpolation property we specialize to the case that $F$ is imaginary quadratic. Let $p$ be a rational prime that is split in $F$. Let $\lambda$ be a Hecke character of infinity type $(1,0)$. Then Katz constructed a $p$-adic meromorphic function
$$L_p(\chi, t, s):\Z_p\times\Z_p\longrightarrow\C_p$$
determined by the following interpolation property: for all $k,j\in\Z$ such that $1\leq j\leq k-1$,
\begin{equation}\label{KatzInt}\frac{L_p(\chi, k, j)}{\Omega_p^{k-1}} =  E_p(\chi, k, j)\frac{L(\chi\lambda^{k-1}, j)}{\Omega_\infty^{k-1}}.\end{equation}
Here $E_p(\chi, k, j)$ is an explicit complex number and $\Omega_p\in\C_p^\times$, $\Omega_\infty\in\C^\times$ are $p$-adic and complex periods that make both sides of (\ref{KatzInt}) algebraic.

In these two cases, $F$ totally real and $F$ imaginary quadratic with $p$ split, $p$-adic Stark conjectures have been made for $L_p(\chi,s)$ and $L_p(\chi, t, s)$, and some progress has been made on these conjectures. When $F$ is totally real and $\chi$ is totally odd Gross (\cite{Gr}) stated a conjecture for the order of vanishing of $L_p(\chi\omega,s)$ at $s = 0$ and the leading term of the Taylor series of $L_p(\chi\omega, s)$ at $s=0$. Progress has been made on the order of vanishing, and recently the formula for the leading term was proved in \cite{DKV} building off of earlier work in \cite{DDP}. When $F$ is totally real and $\chi$ is totally even there is a conjecture for the value $L_p(\chi,1)$ known as the Serre-Solomon-Stark conjecture (\cite{So}, \cite{Ta}). This conjecture is open except in the cases when $F = \Q$, when the formula is due to Leopoldt, and when $\chi$ is trivial, where Colmez has proven a $p$-adic class number formula (\cite{Co}). When $F$ is imaginary quadratic and $p$ is split in $F$, Katz stated and proved a $p$-adic Stark conjecture for the value $L_p(\chi, 1, j)$ known as Katz's $p$-adic Kronecker's 2nd limit formula (\cite{K1} and see Section \ref{KatzpAdicLfunction}).

One of the original motivations for Stark's conjectures is that when the order of vanishing of $L(\chi,s)$ at $s=0$ is exactly one, then the conjectures shed light on Hilbert's 12th problem about explicit class field theory. More precisely, when the order vanishing is exactly one then Stark's conjectures predict the existence of a unit $u\in\OX_K^\times$ such that the relation
\begin{equation}\label{Stark1}L'(\psi,0) = -\frac{1}{e}\sum_{\sigma\in\Gal(K/F)}\psi(\sigma)\log|\sigma(u)|\end{equation}
holds for all characters of the Galois group $\Gal(K/F)$ and such that $K(u^{1/e})$ is an abelian extension of $F$. Here $e$ is the number of roots of unity in $K$ and the absolute value is a particular absolute value on $K$. When $F$ is real quadratic, $\ord_{s = 0}(L(\chi,s)) = 1$ if and only if $\chi$ is mixed signature. In this case, we choose the absolute value on $K$ to correspond to one of the real places of $K$. Then by varying $\psi$ and exponentiating (\ref{Stark1}) one can solve for the unit $u$ from the $L$-values $L'(\psi,0)$. In this way, Stark's conjectures give a way to construct units in abelian extensions of $F$. In Section \ref{RankOneAb}, we review the rank one abelian Stark conjecture when $F$ is a quadratic field.

The goal of this article is to define a $p$-adic $L$-function and state a $p$-adic Stark conjecture in the setting when $F$ is a quadratic field and $\ord_{s=0}(L(\chi,s)) = 1$ (the rank one setting). This is the case when $\chi$ is any nontrivial character if $F$ is imaginary quadratic, and when $\chi$ is a mixed signature character when $F$ is real quadratic. When $F$ is imaginary quadratic and $p$ is split in $F$ our $p$-adic $L$-function is related to Katz's. In the cases when $F$ is imaginary quadratic and $p$ is inert, as well as when $F$ is real quadratic and $\chi$ is mixed signature, our $p$-adic $L$-function is new. One of the main issues with defining the $p$-adic $L$-function for $\chi$ when $F$ is quadratic and $\ord_{s=0}(L(\chi,s)) = 1$ is that the complex $L$-function $L(\chi,s)$ has no critical values. Therefore the $p$-adic $L$-function of $\chi$ will not interpolate any of the special values of $L(\chi,s)$. In order to define the $p$-adic $L$-function in lieu of the fact that $L(\chi,s)$ has no critical values we $p$-adically deform $\chi$ into a family of $p$-adic representations where complex $L$-functions in the family do have critical values to interpolate.

We now explain in more detail our definition, conjectures, and results. Let
$$\rho = \Ind_{G_F}^{G_\Q}:G_{\Q}\longrightarrow \GL_2(\C)$$
be the induction of $\chi$ from $G_F$ to $G_\Q$. Then the $q$-expansion
$$f = \sum_{\ga\subset\OX_F}\chi(\ga)q^{N\ga}$$
is the $q$-expansion of a weight one modular form and $\rho$ is the representation associated to $f$. The modular form $f$ has character $\varepsilon = \det\rho$ and level $N = |d_F|N_{F/\Q}\mathfrak{m}$ where $d_F$ is the discriminant of $F$ and $\mathfrak{m}$ is the conductor of $\chi$. Let
$$x^2 - a_p(f)x + \varepsilon(p) = (x - \alpha)(x - \beta)$$
be the Hecke polynomial of $f$ at $p$. Then $\alpha$ and $\beta$ are roots of unity, so $f$ has two (possibly equal) ordinary $p$-stabilizations. Let $f_\alpha(z) = f(z) - \beta f(pz)$ be a $p$-stabilization of $f$. Under the assumption that $\alpha\not=\beta$, Bella\"{i}che and Dmitrov (\cite{BD}) have shown that the eigencurve is smooth at the point corresponding to $f_\alpha$. We will use Bella\"{i}che and Dmitrov's result, so we assume $\alpha\not=\beta$ and let $V$ be a neighborhood of $f_\alpha$ on the eigencurve such that the weight map is \'{e}tale at all points of $V$ except perhaps $f_\alpha$. Let $\We$ be weight space. Using the constructions of $\cite{Be1}$ there exists a two-variable $p$-adic rigid analytic function
$$L_p(f_\alpha,z,\sigma):V\times\We\longrightarrow \C_p$$
such that for all classical points $y\in V$, all finite order characters $\psi\in \We(\C_p)$, and all integers $j$, $1\leq j\leq k-1$ where $k$ is the weight of $y$,
\begin{equation}\label{eigenInt}\frac{L_p(f_\alpha, y, \psi^{-1}(\cdot)\langle\cdot\rangle^{j-1})}{\Omega_{p,y}^{sgn(\psi)}} = E_p(f_\alpha,y,\psi,j)\frac{L(g_y,\psi\omega^{j-1},j)}{\Omega_{\infty,y}^{sgn(\psi)}}.\end{equation}
Here $g_y$ is the modular form corresponding to the point $y\in V$, $L(g_y,\psi\omega^{j-1},j)$ is the complex $L$-function of the modular form $g_y$ twisted by the Dirichlet character $\psi\omega^{j-1}$, $E_p(f_\alpha,y,\psi,j)$ is an explicit complex number, and $\Omega_{\infty,y}^\pm,\Omega_{p,y}^\pm$ are $p$-adic and complex periods respectively that make both sides of the equality algebraic. In Section \ref{BackgroundChapter}, we give the background needed in order to define $L_p(f_\alpha,z,\sigma)$.

Conceptually, it makes sense to define the $p$-adic $L$-function of $\chi$ as
$$L_p(\chi,\alpha,s):\Z_p\longrightarrow\C_p$$
$$L_p(\chi,\alpha,s) = L_p(f_\alpha,x,\langle\cdot\rangle^{s-1})$$
where $x\in V$ is the point corresponding to $f_\alpha$. The problem with this definition is that while the function $L_p(f_\alpha,z,\sigma)$ is determined by the above interpolation property, the triple of the function $L_p(f_\alpha,z,\sigma)$, the $p$-adic periods $\Omega_{p,y}^\pm$, and the complex periods $\Omega_{\infty,y}^\pm$ is not canonically defined. The choice of the function $L_p(f_\alpha,z,\sigma)$ may be changed by a $p$-adic analytic function on $V$ for which we would obtain a new function with new $p$-adic and complex periods satisfying the same interpolation formula.  We would like to state a $p$-adic Stark conjecture for the function $L_p(\chi,\alpha,s)$, but because the function is not canonically defined it does not make sense to specify its value at any point with a precise conjecture.

To define a function that does not depend on any choices, we fix two finite order Dirichlet characters $\eta,\psi\in\We(\C_p)$ and define the $p$-adic $L$-function of $\chi$ with the auxiliary characters $\eta$ and $\psi$ as
$$L_p(\chi,\alpha,\psi\omega,\eta\omega,s) = \frac{L_p(f_\alpha,x,\psi^{-1}\omega^{-1}(\cdot)\langle\cdot\rangle^{s-1})}{L_p(f_\alpha,x,\eta^{-1}\omega^{-1}(\cdot)\langle\cdot\rangle^{s-1})}.$$
The function $L_p(\chi,\alpha,\psi\omega,\eta\omega,s)$ does not depend on the choices made to define $L_p(f_\alpha,x,\sigma)$. In Section \ref{DefinitionChapter}, we make the following conjecture for $L_p(\chi,\alpha,\psi\omega,\eta\omega,s)$.

\begin{conj} Let $\eta,\psi\in\We(\C_p)$ be of orders $p^m$ and $p^n$ respectively. Let $M_m$ and $M_n$ be the fixed fields of the kernels of the representations $\rho\otimes\eta$ and $\rho\otimes\psi$ respectively. Let $k_m$ and $k_n$ be the fields obtained by adjoining the values of $\chi$, $\alpha$, and $\zeta_{p^{m+1}}$ and $\zeta_{p^{n+1}}$ respectively to $\Q$. Then there exists units $u_{\chi\eta,\alpha}^*\in k_m\otimes\OX_{M_m}^\times$ and $u_{\chi\psi,\alpha}^*\in k_n\otimes\OX_{M_n}^\times$ such that
$$L_p(\chi,\alpha,\psi\omega,\eta\omega,0) = \frac{(1 - \beta\psi(p))\left(1 - \frac{\psi^{-1}(p)}{\alpha p}\right)\frac{\tau(\psi^{-1})}{p^{n + 1}}}{(1 - \beta\eta(p))\left( 1 - \frac{\eta^{-1}(p)}{\alpha p}\right)\frac{\tau(\eta^{-1})}{p^{m + 1}}}\frac{\log_p(u_{\chi\psi,\alpha}^*)}{\log_p(u_{\chi\eta,\alpha}^*)}$$
where $\tau(\psi^{-1})$ and $\tau(\eta^{-1})$ are the Gauss sums associated to $\psi^{-1}$ and $\eta^{-1}$ respectively.
\end{conj}

In Section \ref{ImagQuadChapter}, we prove our conjecture when $F$ is imaginary quadratic and $p$ is split in $F$ by comparing $L_p(\chi,\alpha,\psi\omega,\eta\omega,s)$ to Katz's $p$-adic $L$-function (Theorem \ref{finishThmKatz}). We also show in Section \ref{ImagQuadChapter} that in this case, it is possible to choose the periods in (\ref{eigenInt}) in such a way as to make $L_p(\chi,\alpha,s)$ canonically defined. It is a goal of future research to explore whether or not this is possible in the other cases.

At the outset of this project, we believed that the units $u_{\chi\eta,\alpha}^*,u_{\chi\psi,\alpha}^*$ would be related to the units appearing in (\ref{Stark1}) for the characters $\chi\eta$ and $\chi\psi$ (see \cite{Ferr} for the precise relation we expected). This is the case when $F$ is imaginary quadratic and $p$ is split (see Section \ref{finishProofKatz}). In the other cases it is not clear what the precise relation is or if there is a relation. In Section \ref{NumEvChap}, we give evidence for our conjecture exploring the relation between $u_{\chi\eta,\alpha}^*$ and $u_{\chi\psi,\alpha}^*$ and the units that would appear in (\ref{Stark1}).

\textbf{Acknowledgments}: I thank my adviser Samit Dasgupta for giving my the ideas that form the basis of this material. This article is based off the work I did in my thesis (\cite{Ferr}). I greatly thank Ralph Greenberg for sharing his work with Nike Vatsal (\cite{GV}) with me, and for reading and helping with an earlier version of this paper. His input was essential in finishing this project. I  thank Rob Pollack and Rob Harron for letting me use their SAGE code for computing overconvergent modular symbols. I also thank the anonymous referee whose comments helped improve the quality of this paper.

\section{The rank one abelian Stark conjecture for quadratic fields}\label{RankOneAb}

In this section we state the rank one abelian Stark conjecture for quadratic fields, and introduce notation that will be used in later sections. Let $F$ be a quadratic extension of $\Q$ and let $K$ be a nontrivial finite abelian extension of $F$. If $F$ is real quadratic assume that one infinite place of $F$ stays real in $K$ and the other becomes complex.

Let $S$ be a finite set of places of $F$ that contains the infinite places and the places that ramify in $K$. Assume that $|S|\geq 2$. Let $S_K$ denote the places of $K$ above those in $S$. Let $v$ denote an infinite place of $K$ such that $v(K)\subset\R$ if $F$ is real quadratic. We also let $v$ denote the infinite place of $F$ that is $v\rvert_F$, so $v\in S$. Let $U_{v,S}$ denote the subgroup of $K^\times$ defined by
$$U_{v, S} = \begin{cases}\{u\in K^\times: |u|_{w'} = 1, \forall w'\text{ such that } w'\rvert_{F}\not= v\rvert_F\} & \text{ if } |S|\geq 3\\
						 \{u\in K^\times: |u|_{w'} = |u|_{w''},\forall w', w''\mid v'\text{ and }|u|_{w} = 1, \forall w\not\in S_K\} & \text{ if } S = \{v, v'\}.\end{cases}$$
Let $e$ denote the number of roots of unity in $K$. Let $L_S(\chi,s)$ be the complex $L$-function associated to $\chi$ with the Euler factors at the primes in $S$ removed.

\begin{conj}\label{ComplexStarkAt0} (Stark \cite{St4} at $s = 0$) There exists $u\in U_{v, S}$ such that for all characters $\chi$ of $\Gal(K/F)$,
$$L'_S(\chi, 0) = -\frac{1}{e}\sum_{\sigma\in\Gal(K/F)}\chi(\sigma)\log|\sigma(u)|_v.$$
\end{conj}

\begin{rem} \begin{enumerate}[nolistsep]

\item Stark conjectured the additional conclusion that $K(u^{1/e})$ is an abelian extension of $F$. For our purposes we will not be considering this part of the conjecture.

\item Stark proved the above conjecture when $F$ is imaginary quadratic (\cite{St4}). The conjecture is open when $F$ is real quadratic.

\item If $|S|\geq 3$, then the element $u\in U_{v,S}$ has its absolute value specified at every infinite place of $K$, so $u$ if it exists is determined up to multiplication by a root of unity. 

\item In the real quadratic case, we can always take $S$ to be the infinite places of $F$ union the places of $F$ that ramify in $K$. In this case, the conjectural $u\in U_{v,S}$ is an actual unit in $\OX_K$. Similarly in the imaginary quadratic case if at least two primes of $F$ ramify in $K$ and we take $S$ to be the infinite place of $F$ union the places of $F$ that ramify in $K$, then the Stark unit $u\in U_{v, S}$ is a unit in $\OX_K$. 

\end{enumerate}
\end{rem}

\begin{defn}\label{StarkUnitDef} Let $K/F$, $S$ and $v$ be as above. An element in $U_{v,S}$ satisfying the above conjecture is called a \textbf{Stark unit for $K/F$} and is denoted $u_K$. If $|S|\geq 3$, then $u_K$ is determined up to multiplication by a root of unity. When $F$ is imaginary quadratic the units $u_K$ will be specified in Section \ref{ComplexStarkUnits}. 
\end{defn}

Now fix a character $\chi$ of $\Gal(K/F)$. We state the rank one abelian Stark conjecture for the one $L$-function $L_S(\chi,s)$. 

We keep the setting and notation as above for $K/F$, $S$, and $v$. Let $\chi$ be a character of $\Gal(K/F)$ such that $\ord_{s = 0}(L_S(\chi,s)) = 1$, and let $k$ be the field obtained by adjoining the values of $\chi$ to $\Q$. Extend $\log|\cdot |_v$ from $U_{v,S}$ to $k\otimes_\Z U_{v,S}$ by $k$-linearity. Let
$$(k\otimes_\Z U_{v,S})^{\chi^{-1}} = \{u\in k\otimes_\Z U_{v,S} : \sigma(u) = \chi^{-1}(\sigma)u, \forall\sigma\in\Gal(K/F)\}$$
be the $\chi^{-1}$ isotypic component of $k\otimes_\Z U_{v,S}$ where $\Gal(K/F)$ acts via its action on $U_{v, S}$. 

\begin{conj}\label{ComplexStarkForChiAt0}(Rational Stark for $\chi$ at $s = 0$). There exists an element $u_{\chi}\in (k\otimes_\Z U_{v, S})^{\chi^{-1}}$ such that
$$L'_S(\chi,0) = \log|u_{\chi}|_v.$$
\end{conj}

\begin{rem}\label{ComplexIntegralImpliesRational}
\begin{enumerate}[nolistsep]\item As it happens with Conjecture \ref{ComplexStarkAt0}, Conjecture \ref{ComplexStarkForChiAt0} is open when $F$ is real quadratic and $\chi$ is mixed signature.

\item Since we are assuming $\ord_{s=0}(L_S(\chi,s)) = 1$, the $k$-dimension of $(k\otimes_\Z U_{v, S})^{\chi^{-1}}$ is one. 

\item Conjecture \ref{ComplexStarkAt0} implies Conjecture \ref{ComplexStarkForChiAt0} by taking
$$u_{\chi} = -\frac{1}{e}\sum_{\sigma\in\Gal(K/F)}\chi(\sigma)\otimes \sigma(u)\in (k\otimes_\Z U_{v,S})^{\chi^{-1}}$$
where $u\in U_{v,S}$ is the unit satisfying Conjecture \ref{ComplexStarkAt0}. 

\end{enumerate}
\end{rem}

\subsection{The imaginary quadratic case}\label{ComplexStarkUnits}

In this section we define the Stark units that exist in the imaginary quadratic case of the rank one abelian Stark conjecture. These units will be used in later sections.

Let $L = \Z\omega_1 + \Z\omega_2\subset\C$ be a lattice in $\C$ with ordered basis so that $\tau = \omega_1/\omega_2$ is in the upper half plane. Define the sigma and delta functions of a complex number $z$ and lattice $L$ to be
$$\sigma(z, L) = z\prod_{\substack{\omega\in L\\ \omega\not=0}}\left(1 - \frac{z}{\omega}\right)e^{\frac{z}{\omega} + \frac{1}{2}\left(\frac{z}{\omega}\right)^2}$$
$$\Delta(L) = \left(\frac{2\pi i}{\omega_2}\right)^{12}e^{2\pi i \tau}\prod_{n = 1}^\infty (1 - e^{2\pi i n\tau})^{24}.$$
Let
$$A(L) = \frac{\omega_1\overline{\omega_2} - \overline{\omega_1}\omega_2}{2\pi i}$$
so $A(L)$ the area of $\C/L$ divided by $\pi$. Further let
$$\eta_1(L) = \omega_1\sum_{n\in\Z}\sum_{m\in\Z, m\not=0}\frac{1}{(m\omega_1 + n\omega_2)^2}$$
and
$$\eta_2(L) = \omega_2\sum_{m\in\Z}\sum_{n\in\Z,n\not=0}\frac{1}{(m\omega_1 + n\omega_2)^2}$$
and define
$$\eta(z,L) = \frac{\omega_1\eta_2 - \omega_2\eta_1}{2\pi i A(L)}\overline{z} + \frac{\overline{\omega_2}\eta_1 - \overline{\omega_1}\eta_2}{2\pi i A(L)}z.$$
Define the fundamental theta function by
$$\theta(z,L) = \Delta(L)\exp(-6\eta(z,L)z)\sigma(z,L)^{12}.$$

We now define Robert's units associated to an integral ideal of an imaginary quadratic field (\cite{Rob}). For the rest of this section we fix the following notation. Let $F$ be any imaginary quadratic field, $\mathfrak{f}$ a non-trivial integral ideal of $F$, $F(\mathfrak{f})$ the ray class field of $F$ of conductor $\mathfrak{f}$, $G_\mathfrak{f} = \Gal(F(\mathfrak{f})/F)$, $f$ the least positive integer in $\mathfrak{f}\cap \Z$, and $w_\mathfrak{f}$ the number of roots of unity in $F$ congruent to $1\bmod\mathfrak{f}$. For a fractional ideal $\ga$ coprime to $\mathfrak{f}$, let $\sigma_{\ga}\in G_{\mathfrak{f}}$ be the image of $\ga$ under the Artin map. Let $S$ be the set consisting of the infinite place of $F$ and the places dividing $\mathfrak{f}$, and let $v$ be the infinite place of $F(\mathfrak{f})$ induced by $\iota_\infty$.

\begin{defn} Define for $\sigma\in G_\mathfrak{f}$, the \textbf{Robert unit} associated to $\sigma$ by
$$E(\sigma) = \theta(1, \mathfrak{f}\mathfrak{c}^{-1})^f$$
where $\sigma_{\mathfrak{c}} = \sigma$. \end{defn}

\begin{pr} (\cite{GiRob}) For all $\sigma\in G_\mathfrak{f}$,
\begin{enumerate}[nolistsep] \item[(i)] $E(\sigma)\in F(\mathfrak{f})$.
\item[(ii)] For all $\sigma'\in G_\mathfrak{f}$, $\sigma'(E(\sigma)) = E(\sigma'\sigma)$.
\item[(iii)] If $\mathfrak{f}$ is divisible by two distinct primes then $E(\sigma)$ is a unit in $F(\mathfrak{f})$. If $\mathfrak{f} = \q^n$ for a prime $\q$ of $F$, then $E(\sigma)$ is a $\q$-unit.\end{enumerate}
\end{pr}

\begin{thm} (Kronecker's second limit formula)  For all characters $\chi$ of $G_\mathfrak{f}$,
$$L'_S(\chi,0) = -\frac{1}{12fw_\mathfrak{f}}\sum_{\sigma\in G_\mathfrak{f}}\log|E(\sigma)|_v.$$
\end{thm}

When Stark stated his conjectures, he recast this theorem using the following lemma.

\begin{lm} (Lemma 9 on page 225 of \cite{St4}) Let $K\subset F(\mathfrak{f})$ be a subfield of $F(\mathfrak{f})$ that is a nontrivial extension of $F$. Let $J\subset G_\mathfrak{f}$ be the subgroup such that $G_\mathfrak{f}/J = \Gal(K/F)$, and define for $\sigma J\in G_\mathfrak{f}/J$
$$E(\sigma J) = \prod_{\sigma'\in \sigma J} E(\sigma') = N_{F(\mathfrak{f})/K}(E(\sigma)).$$
Let $e$ be the number of roots of unity in $K$. Then $E(\sigma J)^e$ is a $12 fw_\mathfrak{f}$ power in $K$.\end{lm}

\begin{defn}\label{ComplexStarkUnitDef} Let $K\subset F(\mathfrak{f})$ be a nontrivial extension of $F$ such that $\Gal(K/F) = G_\mathfrak{f}/J$. Let $e$ be the number of roots of unity in $K$. Define the \textbf{Stark unit} of the extension $K/F$, denoted $u_K$ to be an element of $K$ such that
$$u_K^{12 f w_\mathfrak{f}} = E(J)^e$$
where $E(J) = \prod_{\sigma\in J} E(\sigma)$. Such an element $u_K$ exists by the previous lemma and is unique up to multiplication by a root of unity in $K$. 
\end{defn}

\begin{thm} (\cite{St4} Stark's Conjecture when $F$ is imaginary quadratic) Keeping the notation as in the previous definition, for all characters $\chi$ of $\Gal(K/F)$,
$$L'_S(\chi,0) = -\frac{1}{e}\sum_{\sigma\in\Gal(K/F)}\chi(\sigma)\log|\sigma(u_K)|_v$$
and $K(u_K^{1/e})$ is an abelian extension of $F$.\end{thm}

\section{Background for definition of the $p$-adic $L$-function}\label{BackgroundChapter}

\subsection{Conventions for modular forms and modular symbols}

In this section, we set some notation and conventions that will be fixed throughout for modular forms and modular symbols. We also state some relevant definitions for later reference.

Fix a positive integer $N$ such that $p\nmid N$ and let $\Gamma$ be either $\Gamma_1(N)$ or $\Gamma_1(N)\cap \Gamma_0(p)$. Our Hecke actions are defined via the double coset algebra of $\Gamma$ in $\GL_2(\Q)$. Let $T_\ell$ denote the Hecke operator at $\ell$ for $\ell\nmid Np$. If $\Gamma = \Gamma_1(N)$, let $T_p$ denote the Hecke operator at $p$, while if $\Gamma = \Gamma_1(N)\cap \Gamma_0(p)$, let $U_p$ denote the Hecke operator at $p$. Let $\iota$ denote the operator for the double coset corresponding to $\begin{pmatrix} 1 &0\\ 0 &-1\end{pmatrix}$. For $a\in (\Z/N\Z)^\times$, let $[a]$ denote the diamond operators. Define the Hecke algebra to be the algebra
$$\mathcal{H} = \begin{cases} \Z[T_\ell, \ell\nmid Np, U_p, [a], a\in(\Z/N\Z)^\times] &\text{ if } \Gamma = \Gamma_1(N)\cap\Gamma_0(p)\\
			\Z[T_\ell, \ell\nmid N, [a], a\in(\Z/N\Z)^\times] &\text{ if } \Gamma = \Gamma_1(N).\end{cases}$$
If $\Sigma$ is a subsemigroup of $\GL_2(\Q)$ containing the matrices needed to define $\mathcal{H}$, then we also consider $\mathcal{H}$ as a subalgebra of the double coset algebra of $\Gamma$ in $\Sigma$.

For $k\geq 1$, let $S_k(\Gamma,\overline{\Q})$ denote the space of holomorphic weight $k$ and level $\Gamma$ cusp forms with algebraic $q$-expansions, and let $S_k(N, \varepsilon, \overline{\Q}) \subset S_k(\Gamma_1(N), \overline{\Q})$ be the space of holomorphic cuspforms of level $N$ and nebentypus $\varepsilon$ with algebraic $q$-expansions.  Let
$$S_k(\Gamma,\C_p) = S_k(\Gamma,\overline{\Q})\otimes_{\overline{\Q}}\C_p\text{ and } S_k(\Gamma,\C) = S_k(\Gamma,\overline{\Q})\otimes_{\overline{\Q}}\C$$
and similarly let
$$S_k(N,\varepsilon,\C_p) = S_k(N,\varepsilon,\overline{\Q})\otimes_{\overline{\Q}}\C_p\text{ and } S_k(N, \varepsilon,\C) = S_k(N, \varepsilon, \overline{\Q})\otimes_{\overline{\Q}}\C.$$

Let $\mathcal{F}$ be the set of holomorphic functions $f$ on the upper half plane such that for all $c\in\Prj^1(\Q)$, $\lim_{z\to c}|f(z)| = 0$, where to make sense of the limit, we view $\Prj^1(\Q)$ and the upper half plane as subsets of $\Prj^1(\C)$. For $k\geq 1$, we define the following weight-$k$ action of $\GL_2^+(\Q)$ on $\mathcal{F}$: for $\gamma = \begin{pmatrix} a & b\\ c & d\end{pmatrix}\in\GL_2^+(\Q)$, $f\in\mathcal{F}$,
$$f\rvert_{\gamma,k}(z) = (cz + d)^kf\left(\frac{az + b}{cz + d}\right).$$
The space $S_k(\Gamma,\C)$, of holomorphic cusp forms of weight $k$ and level $\Gamma$ is the set of invariants of $\Gamma$ with respect to this weight-$k$ action. Let $\Sigma^+ = \GL_2^+(\Q)\cap M_2(\Z)\subset \GL_2(\Q)$. The action of $\Sigma^+$ on $\mathcal{F}$ induces an action of $\mathcal{H}$ on $S_k(\Gamma,\C)$ which leaves the space $S_k(\Gamma,\overline{\Q})$ invariant, defining an action of $\mathcal{H}$ on $S_k(\Gamma,\overline{\Q})$. We extend this action to $S_k(\Gamma,\C_p)$ by linearity. 

For the rest of this article, we adopt the notation that $\Gamma = \Gamma_1(N)$ and $\Gamma_0 = \Gamma_1(N)\cap \Gamma_0(p)$.

\begin{defn} A \textbf{Hecke eigenform} (or just \textbf{eigenform}) of level $N$ and character $\varepsilon$ is an element $f\in S_k(N,\varepsilon,\C_p)$ which is an eigenvector for all the elements of $\mathcal{H}$. A \textbf{normalized eigenform} is a Hecke eigenform $f\in S_k(N,\varepsilon,\C_p)$ such that the leading term of the $q$-expansion of $f$ is $1$. If $f$ is a normalized eigenform, then $f\in S_k(N,\varepsilon,\overline{\Q})$ and so we may also view $f$ as an element of $S_k(N,\varepsilon,\C)$. If $f\in S_k(N,\varepsilon,\overline{\Q})$ is a normalized eigenform that is new at level $N$, we call $f$ a \textbf{newform}.
\end{defn}

\begin{defn} Let $f = \ds\sum_{n = 1}^\infty a_nq^n\in S_k(N, \varepsilon, \overline{\Q})$ be a newform. Then the Hecke polynomial of $f$ at $p$ is the polynomial $x^2 - a_px + \varepsilon(p)p^{k-1}$. Let $\alpha$ and $\beta$ be the roots of this polynomial. Define the \textbf{$p$-stabilizations of $f$} to be $f_\alpha(z) : = f(z) - \beta f(pz)$ and $f_\beta(z) : = f(z) - \alpha f(pz)$.\end{defn} 

The $p$-stabilizations $f_\alpha, f_\beta$ are elements of $S_k(\Gamma_0,\overline{\Q})$, and are eigenvectors for the action of $\mathcal{H}$. The $T_\ell$ eigenvalues of $f_\alpha$ (respectively $f_\beta$) are the same as for $f$ when $\ell\not=p$, and the $U_p$-eigenvalue of $f_\alpha$ (respectively $f_\beta$) is $\alpha$ (respectively $\beta$).

\begin{defn} Let $S_k^{ord}(N,\varepsilon,\C_p)$  (respectively $S_k^{ord}(\Gamma_0,\C_p)$) denote the maximal invariant subspace of $S_k(N,\varepsilon,\C_p)$ (respectively $S_k(\Gamma_0,\C_p)$) with respect to the action of $T_p$ (respectively $U_p$) such that the characteristic polynomial of $T_p$ (respectively $U_p$) restricted to this subspace has roots which are $p$-adic units. We call the subspace $S_k^{ord}(N,\varepsilon,\C_p)$ (respectively $S_k^{ord}(\Gamma_0,\C_p)$) the \textbf{ordinary subspace} of $S_k(N,\varepsilon,\C_p)$ (respectively $S_k(\Gamma_0,\C_p)$). A cuspform $f\in S_k(N, \varepsilon,\C_p)$ (respectively  $S_k(\Gamma_0,\C_p)$) is called \textbf{$p$-ordinary} if $f$ is an element of the subspace $S_k^{ord}(N,\varepsilon,\C_p)$ (respectively $S_k^{ord}(\Gamma_0,\C_p)$).
\end{defn}

We remark that if $f\in S_k^{ord}(N,\varepsilon, \C_p)$ is a newform and $k\geq 2$, then there is a unique $p$-ordinary $p$-stabilization of $f$, while if $f\in S_1(N, \varepsilon,\C_p)$ is a weight one newform, then there are two (possibly equal) $p$-ordinary $p$-stabilizations of $f$.




We now introduce modular symbols. Let $\triangle_0 = \Div^0(\Prj^1(\Q))$ be the set of degree zero divisors on $\Prj^1(\Q)$ and view $\triangle_0$ as a $\GL_2(\Q)$-module via the action of linear fractional transformations. Let $V$ be a right $\Gamma$ module. We define a right action of $\Gamma$ on $\Hom(\triangle_0, V)$ via the rule for $\varphi\in \Hom(\triangle_0, V)$, $\gamma\in \Gamma$, and $D\in\triangle_0$:
$$(\varphi\rvert\gamma)(D) = \varphi(\gamma D)\rvert \gamma.$$

\begin{defn} The set of \textbf{$V$-valued modular symbols on $\Gamma$}, denoted $\Symb_\Gamma(V)$, is the set of all $\varphi\in\Hom(\triangle_0, V)$ that are invariant under the action of $\Gamma$. \end{defn}

In the cases we consider, $V$ has an action of a submonoid of $\GL_2(\Q)$ which defines an action of $\mathcal{H}$ on $\Symb_\Gamma(V)$ through a double coset algebra. When $2$ acts invertibly on $V$ and $\iota$ acts on $\Symb_\Gamma(V)$, we get a decomposition of $\Symb_\Gamma(V)$ into the direct sum of the $1$ and $-1$ eigenspaces of $\iota$, denoted $\Symb^+_\Gamma(V),\Symb^-_\Gamma(V)\subset\Symb_{\Gamma}(V)$. If $\varphi\in\Symb_\Gamma(V)$, then we write $\varphi^\pm$ for the projection of $\varphi$ onto $\Symb^\pm_\Gamma(V)$. 

\subsection{Overconvergent modular symbols} \label{weight space}

In this section we introduce overconvergent modular symbols following the notation and conventions of \cite{Be1} and \cite{PS2}.

Let $\We = \Hom_{cont}(\Z_p^\times, \mathbb{G}_m)$ denote weight space as a $\Q_p$-rigid analytic space, and let $\mathcal{R}$ denote the ring of $\Q_p$-rigid analytic functions on $\We$.  Let $\omega\in\We$ be the Teichm\"{u}ller character. For $m$ with $0\leq m\leq p-2$, let $\We_m\subset\We$ denote the subset of $\We$ consisting of characters whose restriction to $\mu_{p-1}\subset\Z_p^\times$ is equal to $\omega^m$.

We give an explicit description of certain admissible open subsets of the $\Q_p$-points of $\We_m$. For any $\kappa\in \We_m(\Q_p)$ and any $r\geq 1$, let $W(\kappa, 1/p^r)$ denote the closed disk of radius $1/p^r$ in $\We_m$ around $\kappa$. Then
$$W(\kappa, 1/p^r)(\C_p) = \{\kappa'\in \We_m(\C_p): |\kappa'(\gamma) - \kappa(\gamma)|\leq 1/p^r\},$$
and $W(\kappa, 1/p^r)$ is an an admissible open subset of $\We_m$. The ring of $\Q_p$-rigid analytic functions on $W(\kappa, 1/p^r)$ is the $\Q_p$-algebra
$$R = \left\{\sum_{n = 0}^\infty a_n(w - (\kappa(\gamma) - 1))^n\in \Q_p[[w - (\kappa(\gamma)-1)]] : |a_np^{rn}|\rightarrow 0\text{ as } n\rightarrow\infty\right\}$$
and $W(\kappa, 1/p^r) = \Sp R\subset\We$. We remark that $R$ is isomorphic to the Tate algebra
$$\Q_p\langle T\rangle = \left\{\sum_{n = 0}^\infty a_nT^n\in \Q_p[[T]] : |a_n|\rightarrow 0 \text{ as } n\rightarrow\infty\right\}$$
by setting $T = (x - (\kappa(\gamma) - 1)/p^r$. The sets $W(\kappa, 1/p^r)$ form a basis of admissible open neighborhoods of $\kappa$ in $\We_m$.

For each $r\in |\C_p^\times| = p^\Q$, let
$$B[\Z_p, r] = \{z\in \C_p : \exists a \in\Z_p, |z - a| \leq r\}.$$
$B[\Z_p,r]$ is the set of $\C_p$-points of the $\Q_p$-rigid analytic space which is the union of the closed unit balls of radius $r$ around each point in $\Z_p$. Let $\A[r]$ be the $\Q_p$-algebra of rigid analytic functions on $B[\Z_p,r]$.  The sup norm on $\A[r]$ makes $\A[r]$ a $\Q_p$-Banach space. Explicitly the norm is given for $f\in\A[r]$ by
$$\|f\|_r = \sup_{z\in B[\Z_p, r]}|f(z)|.$$

Let $\D[r] = \Hom_{\Q_p}(\A[r], \Q_p)$ be the continuous $\Q_p$-dual of $\A[r]$. The space $\D[r]$ is a $\Q_p$-Banach space with norm given by
$$\|\mu\|_r = \sup_{f\in\A[r], f\not=0}\frac{|\mu(f)|}{\|f\|_r}$$
for $\mu\in\D[r]$. For $r_1 > r_2$, restriction of functions gives a map $\A[r_1]\rightarrow \A[r_2]$. This map is injective, has dense image, and is compact. The dual map $\D[r_2]\rightarrow\D[r_1]$ is injective and compact. For any real number $r\geq 0$ define
$$\A^\dagger[r] = \varinjlim_{s>r}\A[s]\text{ and } \D^\dagger[r] = \varprojlim_{s>r}\D[s].$$
We give $\A^\dagger[r]$ the inductive limit topology and $\D^\dagger[r]$ the projective limit topology. For the remainder of this article, we write $\A = \A^\dagger[0]$ and $\D = \D^\dagger[0]$. We remark that $\D$ is the continuous $\Q_p$-linear dual to $\A$, and that $\A$ may be identified with the set of locally analytic functions on $\Z_p$ and $\D$ the set of locally analytic distributions. 

Given $\mu\in\D$, via integration $\mu$ determines a $\Q_p$-rigid analytic function on $\We$, which we call the $p$-adic Mellin transform of $\mu$. We denote the map corresponding to the $p$-adic Mellin transform
$$\mathcal{L}:\D\longrightarrow \mathcal{R}.$$
For $\mu\in\D$ and $\chi\in\We(\C_p)$, we use the integration symbol for $\mu$ evaluated at $\chi$:
$$\mathcal{L}(\mu)(\chi)=\int_{\Z_p^\times}\chi(z)d\mu(z).$$

We now define overconvergent modular symbols. Let
$$\Sigma_0(p) = \left\{\begin{pmatrix}a & b\\ c & d\end{pmatrix}\in M_2(\Z_p) : p\nmid a, p\mid c \text{ and } ad - bc\not=0\right\}.$$
For any integer $k\in\Z$, we define a weight $k$ action of $\Sigma_0(p)$ on $\A[r]$ for $r<p$ as follows. For $\gamma = \begin{pmatrix} a & b\\ c & d\end{pmatrix}\in \Sigma_0(p)$, $f\in \A[r]$, let
$$(\gamma\cdot_k f)(z) = (a + cz)^kf\left(\frac{b + dz}{a + cz}\right).$$
This induces an action of $\Sigma_0(p)$ on $\D[r]$ on the right via
$$(\mu\rvert_k\gamma)(f) = \mu(\gamma\cdot_kf)$$
for $\mu\in\D[r]$. These actions induce actions of $\Sigma_0(p)$ on $\A$ and $\D$. When we consider $\A$ and $\D$ with their weight $k$ actions, we write $k$ in the subscript, $\A_k$, $\D_k$. The spaces of modular symbols of interest are $\Symb_{\Gamma_0}(\D_k)$. These space are Hecke modules via the action of $\Sigma_0(p)$ on $\D_k$.

\begin{defn} Let $k\in\Z$. The space of \textbf{overconvergent modular symbols of weight $k$} is defined to be $\Symb_{\Gamma_0}(\D_k)$.\end{defn}

\begin{defn} Let $\varphi\in \Symb_{\Gamma_0}(\D_k)$ be an overconvergent modular symbol of weight $k$. We define the $p$-adic $L$-function of $\varphi$ by composing the following two maps: first evaluation at $\{0\} - \{\infty\}$, and then the map $\mathcal{L}$ from before. The composition is called the \textbf{Mellin transform} of $\varphi$ and denoted by $\Lambda_k$:
$$\Lambda_k:\Symb_{\Gamma_0}(\D_k)\rightarrow \mathcal{R}.$$
For $\varphi\in\Symb_{\Gamma_0}(\D_k)$ and $\chi\in\We(\C_p)$,
$$\Lambda_k(\varphi)(\chi) = \int_{\Z_p^\times}\chi(z)d(\varphi(\{0\} - \{\infty\}))(z).$$
By definition, $\Lambda_k$ is a $\Q_p$-linear map. \end{defn}

\subsection{The $p$-adic $L$-function of an ordinary weight $k \geq 2$ modular form}\label{pAdicLfunctionModForm}

In this section, we review how to use classical and overconvergent modular symbols to define the $p$-adic $L$-function of  a weight $k+2\geq 2$ $p$-ordinary newform. 

Let $R$ be a $\Q$-algebra, and for $k\in\Z_{\geq 0}$, let $V_k(R) = \Sym^k(R^2)$ be the $R$-module of homogeneous polynomials of degree $k$ in two variables $X$ and $Y$ with coefficients in $R$. Define an action of $\GL_2(R)$ on $V_k(R)$ as follows: for $\gamma = \begin{pmatrix} a & b\\ c & d\end{pmatrix}\in \GL_2(R)$ and $P\in V_k(R)$, define
$$(P\rvert \gamma)(X,Y) = P((X, Y)\gamma*) = P(dX - cY, -bX + aY)$$
where $\gamma* = \begin{pmatrix} d &-b\\ -c & a\end{pmatrix}$. Since $R$ is a $\Q$-algebra, the space of modular symbols, $\Symb_{\Gamma_0}(V_k(R))$, obtains an action of $\GL_2(\Q)$, which determines a Hecke action of $\mathcal{H}$.

Let $g\in S_{k + 2}(\Gamma_0,\C)$. Define the standard modular symbol associated to $g$, denoted $\psi_g$, to be the function
$$\psi_{g}:\triangle_0\longrightarrow V_k(\C)$$
$$\psi_{g}(\{b\} - \{a\}) = 2\pi i \int_a^bg(z)(zX + Y)^{k}dz.$$
It follows that $\psi_g\in \Symb_{\Gamma_0}(V_k(\C))$ and the map
$$\begin{array}{rcl} S_{k+2}(\Gamma_0,\C) &\longrightarrow &\Symb_{\Gamma_0}(V_k(\C))\\
					g &\longmapsto &\psi_g\end{array}$$
is Hecke equivariant.

Let $f\in S_{k + 2}^{ord}(N,\varepsilon,\C)$ be a $p$-ordinary newform, and let $f_\alpha\in S_{k + 2}^{ord}(\Gamma_0,\C)$ be its $p$-ordinary $p$-stabilization. Shimura (\cite{Sh2}) showed that there exist complex periods $\Omega_{f_\alpha}^{\pm}\in\C^\times$ such that $\psi_{f_\alpha}^\pm/\Omega_{f_\alpha}^{\pm}\in\Symb_{\Gamma_0}^\pm(V_k(\overline{\Q}))$, and that the Hecke eigenspaces in $\Symb^\pm_{\Gamma_0}(V_k(\overline{\Q}))$ with the same eigenvalues as $f_\alpha$ are one-dimensional over $\overline{\Q}$. 

The algebraicity result of Shimura allows one to view the modular symbol associated to $f_\alpha$ $p$-adically in order to define the $p$-adic $L$-function of $f_\alpha$. Let $\varphi_{f_\alpha}^\pm = \psi_{f_\alpha}^\pm/\Omega_{f_\alpha}^\pm\in \Symb_{\Gamma_0}(V_k(\overline{\Q}))$ for some choice of complex periods $\Omega_{f_\alpha}^\pm$. (Each choice of period is determined up to a scalar in $\overline{\Q}^\times$.) Via $\iota_p$, view $\varphi_{f_\alpha}^\pm$ as an element of $\Symb_{\Gamma_0}(V_k(\C_p))$.

Mazur-Tate and Teitelbaum (\cite{MTT}) proved that the function $\mu^\pm_{f_\alpha}$ defined by the rule 
$$\mu^\pm_{f_\alpha}(a + p^m\Z_p) = \alpha^{-m}\varphi_{f_\alpha}^\pm\left(\left\{\frac{a}{p^m}\right\} - \{\infty\}\right)\rvert_{X = 0, Y = 1}$$
is a $\C_p$ valued measure on $\Z_p^\times$. Given a finite order character $\psi\in\We(\C_p)$, we then define the $p$-adic $L$-function of $f_\alpha$ twisted by $\psi$ to be the analytic function of $s\in\Z_p$ given by the formula
$$L_p(f_\alpha, \psi, s) = \int_{\Z_p^\times} \psi^{-1}(t)\langle t\rangle^{s - 1}d\mu_{f_\alpha}^{sgn(\psi)}(t).$$

We record here the interpolation property of $L_p(f_\alpha, \psi, s)$ for future reference.

\begin{thm}(\cite{MTT}) Let $f_\alpha$ be the ordinary $p$-stabilization of a $p$-ordinary newform of level $N$ and weight $k+2\geq 2$. Let $\psi\in\We(\C_p)$ be a finite order character of conductor $p^m$. Then $L_p(f_\alpha,\psi,s)$ is a $p$-adic analytic function on $\Z_p$ with the interpolation property that for all integers $j$ with $0<j<k+2$,
$$L_p(f_\alpha, \psi, j) = \frac{1}{\alpha^m}\left(1 - \frac{\psi^{-1}\omega^{1-j}(p)}{\alpha p^{1-j}}\right)\frac{p^{m(j-1)}(j-1)!\tau(\psi^{-1}\omega^{1-j})}{(2\pi i)^{j-1}}\frac{L(f_\alpha, \psi\omega^{j-1}, j)}{\Omega_{f_\alpha}^{sgn(\psi)}}.$$
Here $\tau(\psi^{-1}\omega^{1-j})$ is the Gauss sum associated to $\psi^{-1}\omega^{1-j}$.

\end{thm}

\begin{rem} 
If $f$ is a non-ordinary newform with Hecke polynomial
$$x^2 - a_p(f)x + \varepsilon(p)p^{k+1} = (x - \alpha)(x - \beta)$$
then one may define the $p$-adic $L$-function of either $p$-stabilization $f_\alpha$ or $f_\beta$ of $f$ in the same way as above but a little more care is needed because the distribution $\mu_{f_\alpha}$ (or $\mu_{f_\beta}$) is not a measure. For the critical $p$-stabilization $f_\beta$ when $f$ is $p$-ordinary, even more care is needed. See \cite{PS1} and \cite{Be1} for more information about these cases.
\end{rem}

When $f$ is a weight one modular form there is no modular symbol associated to $f$ and so the above constructions do not work. It is for this reason that we consider overconvergent modular symbols of arbitrary integer weight $k\in\Z$. We now explain the connection between overconvergent modular symbols of weight $k\in\Z_{\geq0}$ and the modular symbols just considered. 

Let $k\in\Z_{\geq 0}$ and define the map
$$\rho_k:\D_k\longrightarrow V_k(\Q_p)$$
$$\rho_k(\mu) = \int_{\Z_p}(Y - zX)^kd\mu(z).$$
The integration in the definition of $\rho_k$ takes place coefficient by coefficient. The map $\rho_k$ is $\Sigma_0(p)$-equivarient, so induces a Hecke equivariant map
$$\rho_k^*:\Symb_{\Gamma_0}(\D_k)\longrightarrow \Symb_{\Gamma_0}(V_k(\Q_p)).$$
Let $\Symb_{\Gamma_0}(\D_k)^{<k + 1}$ (respectively $\Symb_{\Gamma_0}(V_k(\Q_p))^{<k+1}$) denote the subspace of $\Symb_{\Gamma_0}(\D_k)$ (respectively $\Symb_{\Gamma_0}(V_k(\Q_p))$) which is the span of the set of eigenvectors of $U_p$ with eigenvalue that has $p$-adic valuation less than $k+1$. 

\begin{thm}\label{StevensControl} (Stevens' control theorem \cite{PS2}) For $k\in\Z_{\geq 0}$ the map
$$\rho_k^*: \Symb_{\Gamma_0}(\D_k)^{<k+1}\longrightarrow \Symb_{\Gamma_0}(V_k(\Q_p))^{<k+1}$$
is an isomorphism of Hecke modules.
\end{thm}

\begin{rem}  By Theorem \ref{StevensControl}, there exists unique $\widetilde{\varphi}_{f_\alpha}^\pm\in \Symb_{\Gamma_0}(\D_k)$ such that $\rho_k^*(\widetilde{\varphi}_f^\pm) = \varphi_f^\pm$, and we have the compatibility of $p$-adic $L$-function:
$$\Lambda_k(\widetilde{\varphi}^{sgn(\psi)})(\psi^{-1}\langle\cdot\rangle^{s-1}) = L_p(f_\alpha,\psi, s)$$
for all $s\in\Z_p$ and finite order characters $\psi\in\We(\C_p)$.
\end{rem}

\subsection{Families of overconvergent modular symbols}

In this section we introduce families of overconvergent modular symbols, and construct the ordinary locus of the eigencurve over certain open subsets of weight space following \cite{Be1} and \cite{DuEtAll}. 

Embed $\Z$ into $\We(\Q_p)$ by identifying $k\in\Z$ with the map $x\mapsto x^k$ in $\We(\Q_p)$. Let $k'\in\Z$ and let $W = W(k', 1/p^d)\subset\We$ for some $d\in\Z_{\geq 1}$. We will construct the ordinary locus of the eigencurve over the open set $W$.

Let $\gamma$ be a topological generator of $1 + p\Z_p$. Let $R$ be the ring of $\Q_p$-affinoid functions on $W$ in the variable $w$ induced by the choice $\gamma$. Given $k\in\Z\cap W(\Q_p)$ and $F(w)\in R$, define the evaluation at $k$ map:
$$ev_k:R\longrightarrow \C_p$$
$$ev_k(F) = F(\gamma^k-1).$$
For $r\in p^{\Q}$, define
$$\A[r](R) : = \A[r]\widehat{\otimes}_{\Q_p}R$$
The evaluation maps induce maps
$$ev_k:\A[r](R)\longrightarrow\A[r]$$
for all $r$. 

We define an action of $\Sigma_0(p)$ on $\A[r](R)$ that is compatible with the evaluation maps and the action defined in the previous section. Note that $\A[1]$ is the Tate algebra. Let $z$ be the variable for $\A[1]$, so $\A[1] = \Q_p\langle z\rangle$. For $r\leq 1$ the inclusion $\A[1]\hookrightarrow \A[r]$ induces an inclusion $\A[1](R)\hookrightarrow\A[r](R)$. Let $\log_\gamma(u) = \frac{\log_p(u)}{\log_p(\gamma)}$. Define for $\alpha = \begin{pmatrix} a & b\\ c & d\end{pmatrix}\in\Sigma_0(p)$, $0\leq m\leq p-1$,
$$K_{\alpha, m}(z,w) = \omega(a)^m\sum_{n = 0}^\infty\binom{\log_\gamma(\frac{a + cz}{\omega(a)})}{n}w^n = \omega(a)^m(1 + w)^{\log_\gamma(\frac{a + cz}{\omega(a)})}\in\Q_p[[z,w]].$$
By Lemma 2.6 of \cite{DuEtAll}, $K_{\alpha, m}(z,w)\in \A[1](R)$. By construction, for all $k\in\Z\cap W(\Q_p)$,
$$ev_k(K_{\alpha, m}(z,w)) = (a + cz)^k.$$
Let $r$ be such that $r<p^{\frac{p-2}{p-1}}$. View $K_{\alpha,m}$ as an element of $\A[r](R)$ via the inclusion $\A[1](R)\subset\A[r](R)$. Define, using the ring structure of $\A[r](R)$, for $f\in\A[r]$, $F\in R$, $\alpha = \begin{pmatrix} a & b\\ c & d\end{pmatrix}\in\Sigma_0(p)$,
$$\alpha\cdot(f(z) \otimes F(w)) = K_{\alpha,m}(z,w)(f\left(\frac{b + dz}{a + cz}\right)\otimes F(w))$$
on simple tensors and extend this to an action on $\A[r](R)$ by linearity. 

Define
$$\D[r](R) : = \D[r]\widehat{\otimes}_{\Q_p} R,$$
and define an action of $\Sigma_0(p)$ on $\D[r](R)$ as follows: $\D[r]$ is an $\A[r]$-module via
$$(g\cdot\mu)(f) = \mu(gf)$$
where $f,g\in\A[r]$, $\mu\in\D[r]$. Then $\D[r](R)$ is an $\A[r](R)$-module. Define for $\mu\otimes F\in\D[r](R)$, $\alpha = \begin{pmatrix} a & b\\ c & d\end{pmatrix}\in\Sigma_0(p)$,
$$(\mu\otimes F)\rvert\alpha = K_{\alpha, m}(z,w)(\mu\rvert_0\alpha\otimes F)$$
where $\mu\rvert_0\alpha$ is the weight $0$ action on $\D[r]$.

Now let $\D(R) = \varprojlim_{r>0}\D[r](R)$. The actions of $\Sigma_0(p)$ on $\D[r](R)$ induce an action on $\D(R)$. By Lemma 3.2 of \cite{Be1}, there's a natural isomorphism
$$\D\widehat{\otimes}_{\Q_p}R\longrightarrow\D(R).$$
The map induced by evaluation at $k$ from $\D(R)$ to $\D_k$, will be called specialization to weight $k$ and denoted by $sp_k$:
$$\begin{array}{rcl}sp_k :\D(R) &\longrightarrow &\D_k\\
				\mu\otimes F &\longmapsto &ev_k(F)\mu\end{array}.$$
The map $sp_k$ is $\Sigma_0(p)$-equivariant and induces a Hecke equivariant specialization map which we denote by the same name
$$sp_k:\Symb_{\Gamma_0}(\D(R))\longrightarrow\Symb_{\Gamma_0}(\D_k).$$

To end this section, we summarize some results of Bella\"{i}che (\cite{Be1}) about the relation between $\Symb_{\Gamma_0}(\D(R))$ and $\Symb_{\Gamma_0}(\D_k)$ as Hecke modules. 

\begin{defn} Fix a weight $k\in\Z\cap W(\Q_p)$. Let $\Symb_{\Gamma_0}^\pm(\D_k)^o\subset\Symb_{\Gamma_0}^\pm(\D_k)$ (respectively $\Symb_{\Gamma_0}^\pm(\D(R))^o\subset\Symb_{\Gamma_0}^\pm(\D(R))$) be the subspace where $U_p$ acts with slope bounded by $0$ in the sense of \cite{Be1} Section 3.2.4. Let $\mathbb{T}_k^\pm$ (respectively $\mathbb{T}_{W}^\pm$) be the $\Q_p$-subalgebra of $\End_{\Q_p}(\Symb_{\Gamma_0}^\pm(\D_k)^o)$ (respectively the $R$-subalgebra of $\End_{R}(\Symb_{\Gamma_0}^\pm(\D(R))^o)$) generated by the image of $\mathcal{H}$. We call $\Symb_{\Gamma_0}^\pm(\D_k)^o$ (respectively $\Symb_{\Gamma_0}^\pm(\D(R))^o$) the \textbf{ordinary subspace} of $\Symb_{\Gamma_0}^\pm(\D_k)$ (respectively $\Symb_{\Gamma_0}^\pm(\D(R))$).\end{defn}

We have (\cite{Be1} Section 3.2.4) that $\Symb_{\Gamma_0}^\pm(\D(R))^o$ is a finite projective $R$-module. Then since $\T_W^\pm$ is a finite $R$-algebra, $\T_W^\pm$ is an affinoid algebra. Furthermore, $\T_W^\pm$ is torsion-free as an $R$-module and since $R$ is a principal ideal domain, $\T_W^\pm$ is flat.

\begin{thm} (Bella\"{i}che's specialization theorem (Corollary 3.12 in \cite{Be1})) Let $k\in\Z\cap W(\Q_p)$. The specialization map restricted to the ordinary subspaces
\begin{equation} sp_k:\Symb_{\Gamma_0}^\pm(\D(R))^o\longrightarrow \Symb_{\Gamma_0}^\pm(\D_k)^o\end{equation}
is surjective.\end{thm}

Since $sp_k$ is an $\mathcal{H}$-equivariant surjective map, it induces an $\mathcal{H}$-equivariant map $sp_k:\T_{W}^\pm\rightarrow \T_k^\pm$, which we use in the following definition.

\begin{defn} Let $x:\mathbb{T}_k^\pm\rightarrow\C_p$ be a $\Q_p$-algebra homomorphism. The homomorphism $x$ corresponds to a system of $\mathcal{H}$-eigenvalues appearing in $\Symb_{\Gamma_0}^\pm(\D_k)^o$. Let $\Symb_{\Gamma_0}^\pm(\D_k)_{(x)}$ denote the corresponding generalized eigenspace and let $\Symb_{\Gamma_0}^\pm(\D_k)[x]$ denote the eigenspace. 
\begin{enumerate}[nolistsep] 

\item Let $(\T_{k}^\pm)_{(x)}$ be the localization of $\T_k^\pm\otimes_{\Q_p}\C_p$ at the kernel of $x$. We have that
$$\Symb_{\Gamma_0}^\pm(\D_k)_{(x)} = \Symb_{\Gamma_0}^\pm(\D_k)^o\otimes_{\T_{k}^\pm}(\T_k^\pm)_{(x)}.$$

\item Through the specialization map, $x$ induces a $\Q_p$-algebra homomorphism which we also denote by $x$:
$$x = x\circ sp_k:\T_{W}^\pm\longrightarrow \C_p.$$
Let $(\T_{W}^\pm)_{(x)}$ be the rigid analytic localization of $\T_{W}^\pm\otimes_{\Q_p}\C_p$ at the kernel of $x\circ sp_k$, and let
$$\Symb_{\Gamma_0}^\pm(\D(R))_{(x)} = \Symb_{\Gamma_0}^\pm(\D(R))^o\otimes_{\T_{W}^\pm}(\T_{W}^\pm)_{(x)}.$$

\end{enumerate}

Let $R_{(k)}$ be the rigid analytic localization of $R\otimes_{\Q_p}\C_p$ at the kernel of $ev_k$. We can then localize the specialization map $sp_k$ to get a map
$$sp_k:(\T_{W}^\pm)_{(x)}\otimes_{R_{(k)}, k}\C_p\longrightarrow(\T_{k}^\pm)_{(x)}.$$

\end{defn}

In (\cite{Be1}), Bella\"{i}che following Stevens uses these spaces of families of overconvergent modular symbols to construct the eigencurve. Let $\mathcal{C}_W^\pm = \Sp\T_W^\pm$. Then $\mathcal{C}_W^\pm$ is the ordinary locus of the eigencurve above the open set $W$ of weight space. The weight map 
$$\kappa^\pm:\mathcal{C}_W^\pm\longrightarrow W$$
is the map of rigid analytic spaces induced by the $\Q_p$-algebra homomorphism $R\rightarrow \T_W^\pm$. Since $\T_W^\pm$ is a finite, flat $R$-module, the map $\kappa^\pm$ is finite and flat. Given a point $x\in C_W^\pm(\C_p)$, we define the weight of $x$ to be $\kappa^\pm(x)\in W(\C_p)$. For any $\kappa\in W(\C_p)$ we may consider the evaluation at $\kappa$ map
$$ev_\kappa:R\longrightarrow \C_p$$
$$ev_\kappa(F) = F(\kappa(\gamma) - 1).$$
Define $R_{(\kappa)}$ to be the rigid analytic localization of $R$ at the kernel of $ev_\kappa$.

\begin{thm} (\cite{Be1}) \label{BelProp} Let $U = W(k', 1/p^d)$ for some $k',d\in\Z$, $d\geq1$, and let $x\in C_{U}^\pm(\C_p)$ be a smooth point of weight $k'$. Then there exists a neighborhood, $W = W(k', 1/p^r)$ of $k'$ with $r\geq d$, such that the following hold. Let $R$ be the ring of rigid analytic functions on $W$. Let $T$ be the direct factor of $\T_W^\pm$ corresponding to the connected component of $C_{W}^\pm$ that $x$ lies in. (Note that $T$ may be defined over a finite extension of $\Q_p$.) Let $T_{\C_p} = T\widehat{\otimes}\C_p$ and $R_{\C_p} = R\widehat{\otimes}\C_p$.

\begin{enumerate}[nolistsep]

\item The generalized eigenspace $\Symb_{\Gamma_0}^\pm(\D_{k'})_{(x)}$ is free of rank one over the algebra $(\T_{k'}^\pm)_{(x)}$, and the eigenspace $\Symb_{\Gamma_0}^\pm(\D_{k'})[x]$ is dimension one over $\C_p$.

\item For all points $y\in C_{W}^\pm$, except perhaps $x$, the algebra $(\T_W^\pm)_{(y)}$ is \'{e}tale over $R_{(\kappa^\pm(y))}$.

\item There exists $u\in R_{\C_p}$ such that $ev_\kappa(u) = 0$ and $\kappa$ is the only $0$ of $u$ on $W$ and an element $t\in T$ such that $x(t)\not=0$ as well as an isomorphism
$$T_{\C_p}\longrightarrow R_{\C_p}[X]/(X^e - u)$$
sending $t$ to $X$.

\item The $T_{\C_p}$-module $\Symb_{\Gamma_0}^\pm(\D(R))^o\otimes_{\T_W^\pm}T_{\C_p}$ is free of rank one.

\item For any point $y\in C_W^\pm(\C_p)$ of weight $\kappa^\pm(y)\in\Z$, the $\mathcal{H}$-equivariant map
$$\Symb_{\Gamma_0}^\pm(\D(R))^o\otimes_{\T_W^\pm}T_{\C_p}\longrightarrow \Symb_{\Gamma_0}^\pm(\D_{\kappa^\pm(y)})_{(y)}$$
sends any generator of $\Symb_{\Gamma_0}^\pm(\D(R))^o\otimes_{\T_W^\pm}T_{\C_p}$ to a generator of $\Symb_{\Gamma_0}^\pm(\D_{\kappa^\pm(y)})_{(y)}$.

\end{enumerate}

\end{thm}
\begin{proof} This theorem is a combination of results from Section 4 of \cite{Be1}. \end{proof}

\begin{rem} If $k\in\Z_{\geq0}$, then the eigencurve is \'{e}tale over weight space at any weight $k$-point. The point that we are interested in is when $k=-1$, which corresponds to weight one modular forms. At weight $k = -1$ points, the eigencurve may not be \'{e}tale over weight space, and this is the case in which we will use the above proposition.
\end{rem}

\subsection{Two-variable $p$-adic $L$-functions}\label{TwoVariable}

In this section we explain how to use Theorem \ref{BelProp} to construct a two-variable $p$-adic $L$-function. 

Let $W = W(k', 1/p^r) = \Sp R$ for some $k'\in \Z$ and $r\geq 1$. Let $M = \Symb_{\Gamma_0}^\pm(\D(R))^o$. Define the $R$-linear map
$$\Lambda:M\longrightarrow R\widehat{\otimes}_{\Q_p}\mathcal{R}$$
to be the composition of evaluation at $\{0\} - \{\infty\}$ and the map $\mathcal{L}$ from before. By construction, for all $k\in\Z\cap W(\Q_p)$ we have the commutative diagram
\begin{equation}\diagram
{M & R\widehat{\otimes}\mathcal{R}\\
 \Symb_{\Gamma_0}^\pm(\D_k)^o & \mathcal{R}.\\};
 \path[->,font = \scriptsize]
 (m-1-1) edge node[auto] {$\Lambda$} (m-1-2)
 (m-1-1) edge node[auto] {$sp_k$} (m-2-1)
 (m-2-1) edge node[auto] {$\Lambda_k$} (m-2-2)
 (m-1-2) edge node[auto] {$ev_k$} (m-2-2);
\end{tikzpicture}\end{equation}
which shows that for $\Phi\in M$, $\Lambda(\Phi)$ interpolates the functions $\Lambda_k(sp_k(\Phi))$.

We now put ourselves in the situation of Theorem \ref{BelProp}, and we extend scalars to $\C_p$. A $\C_p$ in the subscript means completed tensor product over $\Q_p$ with $\C_p$. Let $x\in C_{U}^\pm(\C_p)$ be a smooth point of weight $k'\in \Z$ for some $U = W(k',1/p^d)$. Let $W = W(k', 1/p^r) = \Sp R$ and $T$ be as in the proposition. Let $\epsilon\in \T_{W,\C_p}^\pm$ be such that $T_{\C_p} = \epsilon \T_{W,\C_p}^\pm$. Then
$$\Symb_{\Gamma_0}^\pm(\D(R))^o\otimes_{\T_W^\pm} T_{\C_p} = \epsilon\Symb_{\Gamma_0}^\pm(\D(R))^o_{\C_p}\subset \Symb_{\Gamma_0}^\pm(\D(R))^o_{\C_p},$$
so we let
$$M = \Symb_{\Gamma_0}^\pm(\D(R))^o\otimes_{\T_W^\pm} T_{\C_p} = \epsilon\Symb_{\Gamma_0}^\pm(\D(R))^o_{\C_p}.$$

We first give a construction of a two-variable $p$-adic $L$-function, that we use when the weight map $\kappa^\pm:\mathcal{C}_W^\pm\rightarrow W$ is \'{e}tale. Assume $\kappa^{\pm}:C_W^\pm\rightarrow W$ is \'{e}tale.

The module $M$ is a rank one $T_{\C_p}$-module, so let $\Phi$ be a generator. Let
$$\Lambda(\Phi,\cdot,\cdot):W\times\We\longrightarrow \C_p$$
be the two-variable rigid analytic function that is the image of $\Phi$ in $R\widehat{\otimes}\mathcal{R}$ under $\Lambda$. By the commutative diagram (3), for all $\sigma\in\We$ and $k\in\Z$, 
$$\Lambda(\Phi, k, \sigma) = \Lambda_k(sp_k(\Phi), \sigma).$$

We now consider the non-\'{e}tale case. In the non-\'{e}tale case, if $y\in \mathcal{C}_W^\pm$ is of weight $k\not=k'$, then $sp_k(\Phi)\in\Symb_{\Gamma_0}^\pm(\D_k)^o$ is not in the eigenspace corresponding to $y$. The construction that follows is due to Bella\"{i}che (\cite{Be1}). Let $N = M\otimes_{R_{\C_p}} T_{\C_p}$ and let $V = \Sp T$. Define
$$\Lambda_T: = \Lambda\otimes Id_{T_{\C_p}}:N\longrightarrow (R_{\C_p}\widehat{\otimes}_{\C_p}\mathcal{R}_{\C_p})\otimes_{R_{\C_p}} T_{\C_p}\cong T_{\C_p}\widehat{\otimes}_{\C_p} \mathcal{R}_{\C_p}.$$
Then for $\Phi\in N$, the function $\Lambda_T(\Phi)\in T_{\C_p}\widehat{\otimes}\mathcal{R}_{\C_p}$ is a two-variable rigid analytic function on $V_{\C_p}\times \We_{\C_p}$. For each $y\in V(\C_p)$ of weight $\kappa\in W(\Q_p)$, define the specialization map
$$sp_y:N\longrightarrow \Symb_{\Gamma_0}(\D_\kappa)_{\C_p}^o$$
as the natural map
$$N \longrightarrow N\otimes_{T_{\C_p},y}\C_p.$$
We view $N\otimes_{T_{\C_p},y}\C_p$ as a subset of $\Symb_{\Gamma_0}^\pm(\D_\kappa)_{\C_p}^o$  via
$$\begin{array}{ll} N\otimes_{T_{\C_p},y}\C_p & = (M\otimes_{R_{\C_p}}T_{\C_p})\otimes_{T_{\C_p},y}\C_p\\
									& = M\otimes_{R_{\C_p},ev_\kappa}\C_p\hookrightarrow \Symb_{\Gamma_0}^\pm(\D_\kappa)_{\C_p}^o.\end{array}$$
By construction $sp_y$ is $\mathcal{H}$-equivariant with respect to the action of $\mathcal{H}$ on the first component of $N$. Furthermore, if $\phi\in N$ and $y\in V(\C_p)$ is of weight $k\in\Z$, then (\cite{Be1} Lemma 4.12)
$$\Lambda_T(\Phi)(y,\sigma) = \Lambda_k(sp_y(\Phi))(\sigma).$$

We recall that we have an element $t\in T_{\C_p}$ and $u\in R_{\C_p}$ and an isomorphism
$$T_{\C_p}\longrightarrow R_{\C_p}[X]/(X^e - u)$$
sending $t$ to $X$. Now let $\phi$ be a generator of $M$ as a $T_{\C_p}$ module, and define
$$\Phi = \sum_{i = 0}^{e-1} t^i\phi\otimes t^{e-1-i}\in N.$$

\begin{pr}\begin{enumerate}[nolistsep] \item Let $T_{\C_p}\otimes_{R_{\C_p}}T_{\C_p}$ act on $N$ with the first factor acting on $M$ and the second factor acting on $T_{\C_p}$. Then
$$(t\otimes 1 - 1\otimes t) \Phi = 0.$$
\item Let $y\in \mathcal{C}_W^\pm(\C_p)$ be a point of weight $\kappa\in W(\Q_p)$. Then 
$$sp_y(\Phi) \in \Symb_{\Gamma_0}^\pm(\D_\kappa)[y].$$
We note that if $y\not=x$, then $\Symb_{\Gamma_0}^\pm(\D_\kappa)[y] = \Symb_{\Gamma_0}^\pm(\D_\kappa)_{(y)}$, while if $y = x$ and the ramification index is $e$, $\Symb_{\Gamma_0}^\pm(\D_\kappa)_{(y)}$ is an $e$-dimensional vector space.
\end{enumerate}\end{pr}

\begin{proof} The first part of the proposition is Lemma 4.13 of \cite{Be1} and the second part is Proposition 4.14 of \cite{Be1}. \end{proof}

Define the two-variable $p$-adic $L$-function to be
$$\Lambda_T(\Phi):V_{\C_p}\times\We_{\C_p}\longrightarrow \C_p.$$
To compare this second construction with the first construction when the ramification index is $1$, we note that $T_{\C_p}\cong R_{\C_p}[X]/(X^e - u) = R_{\C_p}$, so
$$N = M\otimes_{R_{\C_p}}T_{\C_p}\cong M$$
and 
$$\Phi = \sum_{i = 0}^{e-1}t^i\phi\otimes t^{e-1-i} = \phi\otimes 1$$
so the second construction reduces to the first one when $e = 1$.

\section{Definition of $p$-adic $L$-functions and $p$-adic Stark Conjecture}\label{DefinitionChapter}

We begin this section by introducing the objects we are working with and setting notation. Let $F$ be a quadratic field of discriminant $d_F$, and let $\chi:G_F\rightarrow\overline{\Q}^\times$
be a nontrivial ray class character of $F$ that is of mixed signature if $F$ is real quadratic. Let $K$ be the fixed field of the kernel of $\chi$ and let $\mathfrak{f}$ be the conductor of $K/F$. Assume that $\iota_\infty(K)\subset\R$ if $F$ is real quadratic. Let $v$ denote the infinite place of $K$ determined by $\iota_\infty$. Let $\rho = \Ind_{G_F}^{G_\Q}\chi:G_\Q\longrightarrow\GL_2(\overline{\Q})$ be the induction of $\chi$ and let $M$ be the fixed field of the kernel of $\rho$. Let $f$ be the weight one modular form associated to $\rho$, so $f$ has level $N = N_{F/\Q}(\mathfrak{f})\cdot|d_F|$ and character $\varepsilon = \det\rho$. The $q$-expansion of $f$ is
$$f = \sum_{\substack{\ga\subset\OX_F\\ (\ga,\mathfrak{f}) = 1}}\chi(\ga)q^{N\ga}$$
and we have that $L(f,s) = L(\chi,s).$ Let 
$$x^2 - a_p(f)x + \varepsilon(p) = (x - \alpha)(x - \beta)$$
be the Hecke polynomial of $f$ at $p$. We note that when $p$ splits in $F$, say $p\OX_F = \p\overline{\p}$, then $\alpha = \chi(\p)$ and $\beta = \chi(\overline{\p})$, and if $p$ is inert, then $\alpha = \sqrt{\chi(p\OX_F)}$ and $\beta = -\sqrt{\chi(p\OX_F)}$. Let $k$ be the field obtained by adjoining the values of $\chi$ along with $\alpha$ and $\beta$ to $\Q$.

We make some assumptions that will be fixed throughout. First we assume that $p\nmid N$, which implies in particular that $p$ does not ramify in $M$. We further assume that $p\nmid[M:\Q]$, and we assume that $\alpha\not=\beta$. With these assumptions, let $f_\alpha(z) = f(z) - \beta f(pz)$ be a fixed $p$-stabilization of $f$.

\subsection{Definition of $p$-adic $L$-function}\label{DefinitionsOfFunctions}

We use the constructions from the previous section to define our $p$-adic $L$-function. In order to do that, we start with the following result of Bella\"{i}che and Dmitrov about the eigencurve at weight one points.

\begin{thm}\label{BelDmit} (\cite{BD}) Let $g$ be a classical weight one newform of level $N$, whose Hecke polynomial at $p$ has distinct roots. Then the eigencurve is smooth at either $p$-stabilization of $g$. Moreover, the eigencurve is smooth but not \'{e}tale over weight space if and only if the representation associated to $g$ is obtained by induction from a mixed signature character of a real quadratic field in which $p$ splits.\end{thm}

By our assumption that $\alpha\not=\beta$ the above theorem implies that the eigencurve is smooth at the point corresponding to $f_\alpha$. We may break our situation into four cases, the cases when $F$ is either imaginary or real quadratic and when $p$ is either inert or split in $F$. In the case when $F$ is real quadratic and $p$ is split the eigencurve is smooth but not \'{e}tale at $f_\alpha$. In the other three cases the eigencurve is \'{e}tale at $f_\alpha$.  We adopt the notation from the previous section except that we base change everything to $\C_p$ and we drop all the $\C_p$ subscripts. Since we will conjecture the value at $s=0$, we consider the minus subspace of modular symbols. Let $T = T_{\C_p}$, $M\subset\Symb_{\Gamma_0}^-(\D(R))_{\C_p}^o$, $N$, and $R = R_{\C_p}$ be as in Section \ref{TwoVariable} where the point of interest $x$ is the point on the eigencurve corresponding to $f_\alpha$. Let $\phi$ be a generator of $M$ as a $T$-module and let
$$\Phi = \sum_{i = 0}^{e-1}t^i\phi\otimes t^{e-1-i}\in N.$$
Let $V = \Sp T$, $\We = \We_{\C_p}$, $W = \Sp(R)$, and let $\Lambda(\Phi) = \Lambda_T(\Phi)$ to make all the notation uniform.

We record the interpolation formulas for our two-variable rigid analytic function
$$\Lambda(\Phi,\cdot,\cdot):V\times\We\longrightarrow\C_p.$$
For each classical point $y\in V$, let $g_y$ be the weight $k\in\Z_{\geq 2}$ $p$-stabilized newform corresponding to $y$. Let $\Omega_{\infty,g_y}\in \C^\times$ be the complex period used to define the $p$-adic $L$-function associated to $g_y$ as in Section \ref{pAdicLfunctionModForm}. Let
$$\varphi_{g_y}\in \Symb_{\Gamma_0}^-(\D_{k-2})_{(y)}$$
be the unique (by Theorem \ref{StevensControl}) modular symbol specializing under $\rho_k^*$ to $$\psi_{g_y}^-/\Omega_{\infty, g_y}\in\Symb_{\Gamma_0}^-(V_{k-2}(\overline{\Q})).$$
Let $\Omega_{p,g_y}\in\C_p^\times$ be the $p$-adic period such that $sp_{y}(\Phi)/\Omega_{p,g_y} = \varphi_{g_y}.$ For each $y$, the period pair $(\Omega_{\infty,g_y},\Omega_{p,g_y})$ viewed as an element of $\C^\times\times\C_p^\times/\overline{\Q}^\times$, where $\overline{\Q}^\times$ is embedded diagonally, does not depend on any choices. 

\begin{pr}\label{BigInterpolation} The two-variable rigid analytic function $\Lambda(\Phi)$ on $V\times\We$ is determined by the following two interpolation properties. First, for all $y\in V$ and all even characters $\sigma\in\We$, $\Lambda(\Phi,y,\sigma) = 0$. Second, for all $y\in V$ corresponding to a $p$-stabilized newform $g_y$ of weight $k\in\Z_{\geq 2}$, and all odd characters $\psi\langle\cdot\rangle^{j-1}\in\We(\C_p)$ where $\psi$ is a finite order character of conductor $p^m$ and $1\leq j\leq k-1$,
\begin{equation}\label{LambdaInterp}\begin{array}{ll}\ds\frac{\Lambda(\Phi, y, \psi\langle\cdot\rangle^{j-1})}{\Omega_{p,g_y}} = &\ds \frac{1}{a_p(g_y)^m}\left(1 - \frac{\psi\omega^{1-j}(p)}{a_p(g_y)p^{1-j}}\right)\frac{p^{m(j-1)}(j-1)!\tau(\psi\omega^{1-j})}{(2\pi i)^{j-1}}\times\\
		&\ds \times\frac{L(g_y,\psi^{-1}\omega^{j-1}, j)}{\Omega_{\infty,g_y}}.\end{array}\end{equation}
This equality takes place in $\overline{\Q}$. Here $\tau(\psi\omega^{1-j})$ is the Gauss sum associated to $\psi\omega^{1-j}$.\end{pr}

\begin{proof} The first interpolation property follows from the fact that $\Phi$ is in the minus subspace for the action of $\iota$. For the second interpolation property, with the way everything is set up, it follows from the fact that
$$\frac{\Lambda(\Phi,y,\sigma)}{\Omega_{p,g_y}} = \frac{\Lambda_k(sp_y(\Phi),\sigma)}{\Omega_{p,g_y}} = L_p(g_y, \psi, j)$$
where $L_p(g_y, \psi, s)$ is defined using that complex period $\Omega_{\infty,g_y}$. \end{proof}
%
%
\begin{rem} At this point, we would like to define the two-variable $p$-adic $L$-function associated to $\chi$ as
\begin{equation}\begin{array}{c} \label{ChipAdic}\ds L_p(\chi,\alpha,\cdot,\cdot):V\times\Z_p\longrightarrow\C_p\\
L_p(\chi,\alpha,y,s) = \Lambda(\Phi,y,\omega^{-1}\langle\cdot\rangle^{s-1}).\end{array}\end{equation}
The $p$-adic $L$-function $L_p(\chi,\alpha,y,s)$ is determined by the above interpolation formula. The first variable is on the eigencurve varying through the $p$-adic family of modular forms passing through $f_\alpha$ and the second variable is the usual cyclotomic variable. To get the one variable $p$-adic $L$-function associated to $\chi$ we would plug the point $x\in V$ that corresponds to $f_\alpha$. It is then natural to make a conjecture for the value $L_p(\chi,\alpha,x,0)$ that is analogous to Conjectures \ref{ComplexStarkAt0} and \ref{ComplexStarkForChiAt0}, replacing the complex logarithm with the $p$-adic logarithm. 

The issue with making the conjecture this way is that the $p$-adic number $L_p(\chi,\alpha,x,0)$ is not canonically defined because we made a choice for $\phi$. The condition on the choice of $\phi$ is that $\phi$ is a generator of $M$ as a $T$-module. If we choose a different generator of $M$ as a $T$-module (changing $\phi$ by an element of $T^\times$) that would change the value $L_p(\chi,\alpha,x,0)$. Therefore as it stands now, we cannot precisely conjecture the value $L_p(\chi,\alpha,x,0)$. 

This issue of the value $L_p(\chi,\alpha,x,0)$ not being canonically defined is a question for further research. One way to approach the problem is to ask whether or not there is a way to canonically choose the periods $(\Omega_{p,g_y},
\Omega_{\infty,g_y})$ so that they determine a two-variable modular symbol $\phi$ which would in turn define the function $L_p(\chi,\alpha,x,s)$ canonically. It is possible to do this in the case when $F$ is imaginary quadratic and $p$ is split in $F$ (see Section \ref{ChoiceOfPeriodsAndComparison}). In this case when $F$ is imaginary quadratic and $p$ is split in $F$ the two-variable $p$-adic $L$-function $L_p(\chi,\alpha,y,s)$ is not canonically defined (it depends on the choice of canonical periods), but the one-variable $p$-adic $L$-function $L_p(\chi,\alpha,x,s)$ is. 

To get around these issues and make a precise conjecture we exploit the fact that in (\ref{LambdaInterp}) the function $\Lambda(\Phi,y,\sigma)$ interpolates the values of the complex $L$-function of $g_y$ twisted by $p$-power conductor Dirichlet characters. Let $\psi\in\We(\C_p)$ be a $p$-power order character. We could then define, generalizing (\ref{ChipAdic}), the $p$-adic $L$-function of $\chi$ twisted by $\psi$ to be
$$L_p(\chi,\alpha,\psi\omega,y,s) = \Lambda(\Phi, y, \psi^{-1}\omega^{-1}\langle\cdot\rangle^{s-1}),$$
and state a $p$-adic Stark conjecture for the value $L_p(\chi,\alpha,\psi\omega,x,0)$. The value $L_p(\chi,\alpha,\psi\omega,x,0)$ is outside the range of interpolation for the function $\Lambda(\Phi,y,\sigma)$, but if it was in the range of interpolation it would be related to $L(f_\alpha,\psi,0)$ at the point $s=0$. We have the relation $L(f,\psi,s) = L(\chi\psi,s)$, and so a conjecture for the value $L_p(\chi,\alpha, \psi\omega,x,0)$ should have the same shape as the conjecture for the value $L'(\chi\psi,0)$ with the complex logarithm replaced with the $p$-adic logarithm. 

Of course, the value $L_p(\chi,\alpha,\psi\omega,x,0)$ has the same issue of not being canonically defined as $L_p(\chi,\alpha,x,0)$, but since we have the flexibility of using finite order characters $\psi\in\We(\C_p)$ we can make a function that is canonically defined. Fix two $p$-power order characters $\eta,\psi\in\We(\C_p)$ and define the function
$$L_p(\chi,\alpha,\eta\omega,\psi\omega,y,s) = \frac{\Lambda(\Phi,y,\eta^{-1}\omega^{-1}\langle\cdot\rangle^{s-1})}{\Lambda(\Phi,y,\psi^{-1}\omega^{-1}\langle\cdot\rangle^{s-1})}.$$
Then $L_p(\chi,\alpha,\eta\omega,\psi\omega,y,s)$ does not depend on the choice of $\phi$ because the indeterminacy of the periods in the interpolation formula (\ref{LambdaInterp}) cancels out. The value $L_p(\chi,\alpha,\eta\omega,\psi\omega,x,0)$ is then canonically defined independent of any choices, and we formulate a conjecture for this value.\end{rem}

\begin{defn} Let $\eta,\psi\in\We(\C_p)$ be two $p$-power order characters. Define the two-variable $p$-adic $L$-function of $\chi$ with the auxiliary characters $\eta$ and $\psi$ as
$$L_p(\chi,\alpha,\eta\omega, \psi\omega,\cdot,\cdot):V\times\Z_p\longrightarrow \C_p\cup\{\infty\}$$
$$L_p(\chi,\alpha,\eta\omega,\psi\omega,y,s) = \frac{\Lambda(\Phi,y,\eta^{-1}\omega^{-1}\langle\cdot\rangle^{s-1})}{\Lambda(\Phi,y,\psi^{-1}\omega^{-1}\langle\cdot\rangle^{s-1})}.$$
The function $L_p(\chi,\alpha,\eta\omega,\psi\omega,y,s)$ does not depend on the choice of $\Phi$.

Define the $p$-adic $L$-function of $\chi$ with the auxiliary characters $\eta$ and $\psi$ as
$$L_p(\chi,\alpha,\eta\omega,\psi\omega,s) = L_p(\chi,\alpha,\eta\omega,\psi\omega,x,s).$$
\end{defn}

\begin{rem}\label{oneVarDef} We may give the definition of $L_p(\chi,\alpha,\eta\omega,\psi\omega,s)$ without making reference to the two-variable $p$-adic $L$-function. The two-variable $p$-adic $L$-function is introduced for two reasons. The first is that it satisfies an interpolation property, while the one-variable function $L_p(\chi,\alpha,\eta\omega,\psi\omega,s)$ does not. The second is that we will use the two-variable $p$-adic $L$-function to prove our conjectures when $F$ is imaginary quadratic and $p$ is split in $F$. 

To define $L_p(\chi,\alpha,\eta\omega,\psi\omega,s)$ without referencing the two-variable $p$-adic $L$-function, we consider the space, $\Symb_{\Gamma_0}^\pm(\D_{-1})^o$, of weight negative one overconvergent modular symbols. Since the eigencurve is smooth at the point $x$ corresponding to $f_\alpha$ the eigenspace $\Symb_{\Gamma_0}^\pm(\D_{-1})[x]$ with the same eigenvalues as $f_\alpha$ is one-dimensional. If $\varphi_{f_\alpha}^\pm$ is a generator of this eigenspace, then $L_p(\chi,\alpha,\eta\omega,\psi\omega,s)$ may be defined as
$$L_p(\chi,\alpha,\eta\omega,\psi\omega,s) = \frac{\Lambda_{-1}(\varphi^{-}_{f_\alpha}, \eta^{-1}\omega^{-1}\langle\cdot\rangle^{s-1})}{\Lambda_{-1}(\varphi^{-}_{f_\alpha},\psi^{-1}\omega^{-1}\langle\cdot\rangle^{s-1})}.$$
Since $\Lambda(\Phi^-,x,\sigma) = \Lambda_{-1}(sp_x(\Phi^-),\sigma)$ and $0\not=sp_x(\Phi^-)$, this definition is the same as the first definition.\end{rem}

\subsection{$p$-adic Conjecture}

For each $n\in\Z_{\geq 0}$, let $\Q_n$ be the $n$th layer of the cyclotomic $\Z_p$ extension of $\Q$, so 
$$\Gal(\Q_n/\Q) = 1 + p\Z_p/1 + p^{n+1}\Z_p\cong\Z/p^n\Z.$$
Let $\Gamma_n = \Gal(\Q_n/\Q)$. Let $M_n$ be the compositum of $M$ and $\Q_n$. Let $\Delta = \Gal(M/\Q)$, and for $n\geq 0$ let $\Delta_n = \Gal(M_n/\Q)$. By our assumption that $p$ does not ramify in $M$ and $p\nmid [M:\Q]$, restriction gives an isomorphism $\Delta_n\cong\Delta\times \Gamma_n$. For any $n\geq 0$, let $v$ denote the infinite place of $M_n$ induced by $\iota_\infty$. Let $U_n = \OX_{M_n}\subset M_n^\times$ if $M_0$ is not the Hilbert class field of $F$ when $F$ is imaginary quadratic. If $F$ is imaginary quadratic and $M_0$ is the Hilbert class field of $F$, let
$$U_n = \{ u\in M_n^\times : |u|_{w'} = |u|_{w''}, \forall w',w''\mid p, |u|_w = 1,\forall w\nmid p, v\}.$$
Let $k_n$ be the field obtained by adjoining the $p^{n + 1}$st roots of unity to $k$. For a character $\eta$ of $\Gamma_n$, let $(\rho\eta)^*$ be the representation $\Ind_{G_F}^{G_\Q}\chi^{-1}\otimes \eta^{-1}$ of $\Delta_n$. Given a $k_n[\Delta_n]$-module $A$, let $A^{(\rho,\eta)*}$ denote the $(\rho\eta)^*$-isotypic component of $A$.

The following is how $\alpha$ is incorporated into our conjectures. It is an idea of Greenberg and Vatsal (\cite{GV}), and is a key aspect to the conjecture. Let $D_p\subset \Delta$ be the decomposition group at $p$ determined by $\iota_p$ and let $\delta_p$ be the geometric Frobenius. For a $k[D_p]$-module $A$, let $A^{\delta_p = \alpha}$ be the subspace where $\delta_p$ acts with eigenvalue $\alpha$. Via the isomorphism $\Delta_n = \Delta\times\Gamma_n$, we view $D_p$ as a subgroup of $\Delta_n$ for any $n$. Then the $\Delta_n$-modules $U_n$ are also $D_p$-modules.

Let $\log_p:\C_p^\times\rightarrow \C_p$ denote Iwasawa's $p$-adic logarithm.  Extend $\log_p$ to $\overline{\Q}\otimes_\Z\C_p^\times$ by $\overline{\Q}$-linearity.

\begin{conj}\label{pAdicStarkAt0} Let $\psi,\eta\in\We(\C_p)$ be of orders $p^n$ and $p^m$ respectively with $m,n\geq 1$. Then there exists units $u_{\chi\psi,\alpha}^*\in (k_n\otimes U_n)^{(\rho\psi)^*, \delta_p = \alpha}$ and $u_{\chi\eta,\alpha}^*\in (k_m\otimes U_m)^{(\rho\eta)^*,\delta_p = \alpha}$ such that
\begin{equation}\label{pAdicStarkAt0Eq}L_p(\chi,\alpha,\psi\omega, \eta\omega, 0 ) = \frac{(1 - \beta\psi(p))\left(1 - \frac{\psi^{-1}(p)}{\alpha p}\right)\frac{\tau(\psi^{-1})}{p^{n + 1}}}{(1 - \beta\eta(p))\left( 1 - \frac{\eta^{-1}(p)}{\alpha p}\right)\frac{\tau(\eta^{-1})}{p^{m + 1}}}\frac{\log_p(u_{\chi\psi,\alpha}^*)}{\log_p(u_{\chi\eta,\alpha}^*)}.\end{equation}
\end{conj}

\begin{rem}\begin{enumerate}[nolistsep]
\item This conjecture should be compared with Conjecture \ref{ComplexStarkForChiAt0}. We are relating the $p$-adic $L$-value $L_p(\chi,\alpha,\psi\omega,\eta\omega,0)$ to the spaces $(k_n\otimes U_n)^{(\psi\rho)^*, \delta_p = \alpha}$ and $(k_m\otimes U_m)^{(\eta\rho)^*, \delta_p = \alpha}$ via the $p$-adic logarithm. The spaces $(k_n\otimes U_n)^{(\psi\rho)^*, \delta_p = \alpha}$ and $(k_m\otimes U_m)^{(\eta\rho)^*, \delta_p = \alpha}$ have $k_n$- and $k_m$-dimension one respectively.

Let $K_n$ be the fixed field of the kernel of $\chi\psi$. At the outset of this project, it was expected that the unit $u_{\chi\eta,\alpha}^*$ would be the projection of the unit $u_{K_n}$ from definition \ref{StarkUnitDef} to the space $(k_n\otimes U_n)^{(\psi\rho)^*,\delta_p = \alpha}$ (\cite{Ferr}). This is the case when $F$ is imaginary and $p$ is split in $F$ (see Section \ref{ImagQuadChapter}), while the numerical evidence suggests that this is not the case when $F$ is imaginary quadratic and $p$ is inert in $F$ (see Sections \ref{Ex2}, \ref{Ex3}). We do think that this is the case when $F$ is real quadratic as we verify in the example in Section \ref{Ex1}, but we do not have enough evidence to conjecture it. 

\item It is also possible to state a conjecture for the $p$-adic value at $s = 1$ (see \cite{Ferr}), and there should be a functional equation relating the two conjectures.

\item In (\cite{GV}), Greenberg and Vatsal define a Selmer group associated to the representation $\rho$ and prove that the characteristic ideal of the Selmer satisfies an interpolation property that is similar to the statement of our conjecture. Proving a main conjecture relating the characteristic ideal of the Selmer group associated to $\rho$ to the analytic $p$-adic $L$-functions defined here would allow one to prove this conjecture using Greenberg and Vatsal's result.

\end{enumerate}
\end{rem}

\section{Proof of the conjecture when $F$ is imaginary quadratic and $p$ splits in $F$}\label{ImagQuadChapter}

\subsection{Katz's $p$-adic $L$-function}\label{KatzpAdicLfunction}

In this section we state relevant facts that are needed about Katz's two variable $p$-adic $L$-function. Let $F$ be an imaginary quadratic field of discriminant $d_F$, and assume $p$ splits in $F$. Let $p$ factor as $p\OX_F = \p\overline{\p}$, where $\p$ is the prime induced by the embedding $\iota_p$. Let $\OX_p = \{x\in\C_p : |x|\leq 1\}$ be the closed unit ball in $\C_p$. Let $\frak{f}$ be an integral ideal of $F$ such that $(\frak{f},p) = 1$. Let $\mathfrak{f}$ factor as $\mathfrak{f} = \prod_{v\mid \mathfrak{f}}\mathfrak{f}_v$. Let $\A_F$ denote the adeles of $F$.

The domain of Katz's $p$-adic $L$-function is the set of all $p$-adic Hecke characters of modulus $\mathfrak{f}$, so we begin by giving our conventions for Hecke characters. Define the subgroups $U_\mathfrak{f}, U_{\mathfrak{f},p}, U_{\mathfrak{f},\infty}\subset\A_F^\times$ as
$$U_\mathfrak{f} = \left\{(x_v)_v\in\A_F^\times: \substack{\ds x_v\equiv 1\bmod \mathfrak{f}_v \text{ if } v\mid \mathfrak{f}\\ \ds x_v\in\OX_{F_v}^\times \text{ if } v\nmid\mathfrak{f}\text{ and is finite}}\right\}$$
$$U_{\mathfrak{f},p} = \left\{(x_v)_v\in\A_F^\times: \substack{\\ \ds x_v\equiv 1\bmod \mathfrak{f}_v \text{ if } v\mid \mathfrak{f}\\ \ds x_v\in\OX_{F_v}^\times \text{ if } v\nmid\mathfrak{f}p\text{ and is finite}\\ \ds x_v = 1 \text{ if } v\mid p}\right\}$$
$$U_{\mathfrak{f},\infty} = \left\{(x_v)_v\in\A_F^\times: \substack{\ds x_v\equiv 1\bmod \mathfrak{f}_v \text{ if } v\mid \mathfrak{f}\\ \ds x_v\in\OX_{F_v}^\times \text{ if } v\nmid\mathfrak{f}\text{ and is finite}\\ \ds x_v = 1\text{ if } v\mid \infty}\right\}.$$
Let $\sigma_1,\sigma_2$ be the two embeddings of $F$ into $\overline{\Q}$. Order $\sigma_1,\sigma_2$ so that $\sigma_1$ is how we view $F$ as a subset of $\overline{\Q}$.

\begin{defn} \begin{enumerate}[nolistsep] \item Let $(a_1,a_2)\in\Z^2$. An \textbf{algebraic Hecke character of $F$ of infinity type $(a,b)$ and modulus $\mathfrak{f}$} is a group homomorphism
$$\chi:\A_F^\times\longrightarrow \overline{\Q}^\times$$
such that the image of $\chi$ is contained in a finite extension of $\Q$, $U_\mathfrak{f}\subset \ker(\chi)$, and for all $x\in F^\times$, $\chi(x) = \sigma_1(x)^{a_1}\sigma_2(x)^{a_2}$. The smallest $\mathfrak{f}$ with respect to divisibility such that $U_{\mathfrak{f}}\subset \ker(\chi)$ is called the \textbf{conductor} of $\chi$. 

If $\chi$ is an algebraic Hecke character of modulus $\mathfrak{f}$ and $\ga$ an ideal of $F$ such that $(\ga,\mathfrak{f}) = 1$ and that factors as $\ds\ga = \prod_{(\p,\ga) = 1} \p^{a_\p}$, then we define $\chi(\ga)$ as
$$\chi(\ga) : = \prod_{(\p,\mathfrak{f}) = 1} \chi(\pi_\p)^{a_\p}$$
where $\pi_\p$ denotes a uniformizer of $F_\p$.

\item A \textbf{$p$-adic Hecke character of $F$} is a continuous group homomorphism
$$\chi:\A_F^\times/F^\times\longrightarrow\C_p^\times.$$
By continuity, there exists an integral ideal $\mathfrak{f}'$ of $F$ such that $(\mathfrak{f}', p) = 1$ and $U_{\mathfrak{f}',p}\subset\ker(\chi)$. Any $\mathfrak{f}'$ for which this is true is called a \textbf{modulus} of $\chi$ and we say that $\chi$ is a $p$-adic Hecke character of modulus $\mathfrak{f}'$. 

\item A \textbf{complex Hecke character of $F$} is a continuous group homomorphism
$$\chi:\A_F^\times/F^\times\longrightarrow \C^\times.$$
By continuity, there exists an integral ideal $\mathfrak{f}'$ of $F$ such that $U_{\mathfrak{f}',\infty}\subset\ker(\chi)$. Any $\mathfrak{f}'$ for which this is true is called a \textbf{modulus} of $\chi$ and we say $\chi$ is a complex Hecke character of modulus $\mathfrak{f}'$. 
\end{enumerate}

If $\chi$ is an algebraic, $p$-adic, or complex Hecke character and $v$ is a place of $F$, then we let $\chi_v$ denote $\chi$ restricted to $F_v^\times\subset \A_F^\times$.
\end{defn}

\begin{rem} In the literature, these notions of Hecke characters go by various names and definitions. We introduce and use the definitions given to avoid confusion. 
\end{rem}

We will also need the following alternative definition of an algebraic Hecke character in terms of ideals. Let $\mathfrak{f}$ be an ideal of $\OX_F$ and let $\alpha\in F^\times$ be an element such that $((\alpha),\mathfrak{m}) = 1$ and say that $\mathfrak{f}$ factors as $\ds\mathfrak{f} = \prod_i \p_i^{f_i}$. Define $\alpha\equiv 1\bmod \mathfrak{f}$ to mean that $\alpha\equiv 1\bmod \p_i^{f_i}$ in $\OX_{F_{\p_i}}$ for all $i$. 

Let $I(\mathfrak{f})$ denote the group of fractional ideals of $F$ that are coprime with $\mathfrak{f}$. Let
$$P_1(\mathfrak{f}) = \{(\alpha)\in I(\mathfrak{f}) : \alpha\in K^\times, \alpha \equiv 1\bmod \mathfrak{f}\}.$$
The second definition of an algebraic Hecke character is, an \textbf{algebraic Hecke character of $F$} of modulus $\mathfrak{f}$ and infinity type $(a_1,a_2)\in\Z^2$ is a group homomorphism $\chi:I(\mathfrak{f})\rightarrow \overline{\Q}^\times$ such that the image of $\chi$ is contained in a finite extension of $\Q$, and for all $\ga\in P_1(\mathfrak{f})$ such that $\ga = (\alpha)$ with $\alpha\equiv 1\bmod\mathfrak{f},$
$$\chi((\alpha)) = \sigma_1(\alpha)^{a_1}\sigma_2(\alpha)^{a_2}.$$
Given an algebraic Hecke character, $\chi$, of modulus $\mathfrak{f}$ and infinity type $(a_1,a_2)$, using the second definition, we get an algebraic Hecke character of the same modulus and infinity type, $\chi_\A$ using the first definition by defining $\chi_\A$ to be the unique group homomorphism $\chi_\A:A_F^\times\longrightarrow \overline{\Q}^\times$ such that:
\begin{enumerate}[nolistsep] \item[(i)] For all primes $\p\in I(\mathfrak{f})$, $\chi_\A\rvert_{\OX_{F_\p}^\times} = 1$ and $\chi_\A(\pi_\p) = \chi(\p)$ for any uniformizer in $F_\p$.
\item[(ii)] For all $x\in F^\times$, $\chi_\A(x) = \sigma_1(x)^{a_1}\sigma_2(x)^{a_2}.$\item[(iii)] $U_\mathfrak{f}\subset\ker(\chi_\A)$.\end{enumerate}

This gives a one-to-one correspondence between algebraic Hecke characters of modulus $\mathfrak{f}$ and infinity type $(a_1,a_2)$ using the first and second definitions.

Given an algebraic Hecke character, $\chi$, of $F$ of infinity type $(a_1,a_2)$ and modulus $\mathfrak{f}$ we obtain $p$-adic and complex Hecke characters $\chi_p$ and $\chi_\infty$ which are defined as follows. Define
$$\chi_p:\A_F^\times/F^\times\longrightarrow \C_p^\times$$
at places $v$ of $F$ not dividing $p$ as $\chi$, so $\chi_p\rvert_{F_v^\times} = \chi\rvert_{F_v^\times}$. At places above $p$ we define $\chi_p$ to be the group homomorphism
$$\chi_p:(F\otimes\Q_p)^\times\longrightarrow\C_p^\times$$
$$\chi_p(\alpha\otimes 1) = \chi(\alpha)/\iota_p(\sigma_1(\alpha)^{a_1}\sigma_2(\alpha)^{a_2}).$$
Since the image of $F^\times$ in $(F\otimes\Q_p)^\times$ is dense this defines $\chi_p$ on $(F\otimes\Q_p)^\times$.  We do something similar for $\chi_\infty$. Define
$$\chi_\infty:\A_F^\times/F^\times\longrightarrow \C^\times$$
at the places $v$ of $F$ not dividing $\infty$ as $\chi$, so $\chi_\infty\rvert_{F_v^\times} = \chi\rvert_{F_v^\times}$. At the place above $\infty$ we define $\chi_\infty$ to be the group homomorphism
$$\chi_\infty:(F\otimes\R)^\times\longrightarrow E_\infty^\times\subset\C^\times$$
$$\chi_\infty(\alpha\otimes1) = \chi(\alpha)/\iota_\infty(\sigma_1(\alpha)^{a_1}\sigma_2(\alpha)^{a_2}).$$
Since the image of $F^\times$ in $(F\otimes\R)^\times$ is dense this defines $\chi_\infty$ on $(F\otimes\R)^\times$.

Given an algebraic Hecke character $\chi$ when we consider $\chi_p$ or $\chi_\infty$, we will drop the subscripts $p$ and $\infty$. It will be clear from context when we are considering $\chi$ as a $p$-adic of complex Hecke character. Furthermore, given a $p$-adic (or complex) Hecke character $\psi$ we may abuse of language and say that $\psi$ is an algebraic Hecke character of infinity type $(a_1, a_2)$ if there exists an algebraic Hecke character $\chi$ of infinity type $(a_1,a_2)$ such that $\psi = \chi_p$ (or $\psi = \chi_\infty$).

Let $\psi$ be an algebraic Hecke character of $F$ of infinity type $(a,b)$ and conductor $\mathfrak{f}'\p^{a_\p}\overline{\p}^{a_{\overline{\p}}}$ where $\mathfrak{f}'$ divides $\mathfrak{f}$. Define the $p$-adic local root number associated to $\psi$ to be the complex number
\begin{equation}\label{pAdicRootNumber} W_p(\psi) = \frac{\psi_\p(\pi_\p^{-a_\p})}{p^{a_\p}}\sum_{u\in(\OX_{F_\p}/\p^{a_\p})^\times}\psi_\p(u)\exp(-2\pi i(Tr_{F_\p/\Q_p}(u/\pi_\p^{a_\p})))\end{equation} 
where $\psi_\p$ denotes $\psi$ restricted to $F_\p^\times$ and $\pi_\p$ is a uniformizer for $F_\p$. Since $F_\p = \Q_p$ we could take $\pi_p = p$.

 Let $G(\mathfrak{f}p^\infty) = \A_F^\times/F^\times U_{\mathfrak{f},p}$, so the space of $p$-adic Hecke characters of $F$ of modulus $\mathfrak{f}$ is
$$\Hom_{cont}(G(\mathfrak{f}p^\infty), \C_p^\times).$$
In \cite{Buzz}, Buzzard explains how to view $\Hom_{cont}(G(\mathfrak{f}p^\infty),\C_p^\times)$ as the $\C_p$-points of a rigid-analytic variety. When we say rigid analytic function in the following theorem it is this rigid analytic structure that we are referring to.

Let $S$ be the set of places containing the infinite places of $F$ and the places of $F$ dividing $\mathfrak{f}$. 

\begin{thm} (\cite{K1}, \cite{deS}) There exists a $p$-adic rigid analytic function 
$$L_p = L_{p,Katz}:\Hom_{cont}(G(\mathfrak{f}p^\infty), \C_p^\times)\longrightarrow \C_p$$
as well as complex and $p$-adic periods $\Omega_\infty\in\C^\times,\Omega_p\in\C_p^\times$ such that for all algebraic Hecke character $\psi$ of $F$ of conductor $\mathfrak{f}'\p^{a_\p}\overline{\p}^{a_{\overline{\p}}}$ where $\mathfrak{f}'$ divides $\mathfrak{f}$ and infinity type $(a,b)$ with $a<0$ and $b\geq 0$, we have
\begin{equation}\label{KatzInterpolation}\frac{L_p(\psi)}{\Omega_p^{b-a}} = \frac{(-a-1)!(2\pi)^{b}}{\sqrt{d_F}^b}W_p(\psi)\left(1 - \frac{\psi^{-1}(\p)}{p}\right)(1 - \psi(\overline{\p}))\frac{L_S(\psi,0)}{\Omega_\infty^{b-a}}.\end{equation}
\end{thm}

\begin{rem}
\begin{enumerate}[nolistsep]
\item Katz originally proved this theorem in \cite{K1} for imaginary quadratic fields and then a similar theorem in \cite{K2} for CM fields. The above statement is taken from \cite{deS} with the correction from \cite{BDP} and with a slight modifications in order to state everything adelically. 

\item The interpolation property (\ref{KatzInterpolation}) uniquely determines Katz's $p$-adic $L$-function.



\end{enumerate}
\end{rem}

We now state Katz's $p$-adic Kronecker's second limit theorem. Let $\zeta_n = \iota_\infty^{-1}(e^{2\pi i/n})\in\overline{\Q}$ for $n\in \Z_{\geq 1}$ be a collection of primitive $n$th roots of unity in $\overline{\Q}$.
\begin{thm}\label{KatzpAdicKron} (\cite{K1}, \cite{deS}) Let $\chi$ be an algebraic Hecke character of conductor $\mathfrak{f}$ and trivial infinity type and let $\psi$ be a Dirichlet character of conductor $p^{n}$. Let $K$ be the fixed field of the kernel of $\chi\psi$ when $\chi\psi$ is viewed as a Galois character via the Artin isomorphism $G(\mathfrak{f}p^\infty) \cong\Gal(F(\mathfrak{f}p^\infty)/F)$. Let $u_K$ be the Stark unit for $K/F$, $G = \Gal(K/F)$, and $e$ be the number of roots of unity in $K$. Then
$$L_p(\chi\psi) = -\frac{1}{e}\frac{\psi(-1)\tau(\psi^{-1})}{\chi(\p^n)p^{n}}\left(1 - \frac{(\chi\psi)^{-1}(\p)}{p}\right)(1 - \chi\psi(\overline{\p}))\sum_{\sigma\in G}\chi\psi(\sigma)\log_p(\sigma(u_K))$$

\end{thm}

\begin{rem} A version of this was proved in Katz's original paper. The formulas for this theorem are taken from \cite{deS} with a minor correction so the $1 - \chi\psi(\overline{\p})$ factor is correct (see \cite{GrFact}). 
\end{rem}

\subsection{Definition of the period pair $(\Omega_\infty,\Omega_p)$}

In this section, we explain how to define the period pair $(\Omega_\infty, \Omega_p)$. The pair $(\Omega_\infty,\Omega_p)$ viewed as an element of $\C^\times\times\C_p^\times/\overline{\Q}^\times$ where $\Q^\times$ is embedded diagonally, is a canonical element associated to $F$. 

Let $K$ be a finite extension of $F$ that contains the Hilbert class field of $F$. Let $\pa$ be the prime of $K$ determined by $\iota_p$. Let $E$ be an elliptic curve with CM by $\OX_F$ defined over $K$ and with good reduction at $\pa$. Let $\omega\in\Omega^1(E/K)$ be an invariant differential of $E$ defined over $K$. Attached to the pair $(E,\omega)$, we let $x$ and $y$ be coordinates on $E$ such that
$$\begin{array}{rcl}\iota:E&\longrightarrow &\Prj^2\\
					P&\longmapsto & (x,y,1)\end{array}$$
is an embedding defined over $K$, which embeds $E$ as the zero set of $y^2 = 4x^3 - g_2x +g_3$ and such that $\iota^*(\frac{dx}{y}) = \omega$. Let $E_\omega$ denote the image of $E$ under $\iota$. Let $E_\omega(\C)\subset \Prj^2(\C)$ denote the complex manifold which consists of the complex points of $E_\omega$. Let $\gamma\in H_1(E_\omega(\C),\Q)$ and define the complex period
$$\Omega_\infty = \frac{1}{2\pi i}\int_\gamma\omega.$$

Let
$$\mathscr{L} = \left\{\frac{1}{2\pi i}\int_\eta\omega: \eta\in H_1(E_\omega(\C),\Z)\right\}$$
be the period lattice of $E_\omega$. We have the complex uniformization
$$\begin{array}{rcl}\Phi:\C/\mathscr{L} &\longrightarrow &E_\omega(\C)\\
						z &\longmapsto &(\mathscr{P}(\mathscr{L},z),\mathscr{P}'(\mathscr{L},z),1)\end{array}$$
where $\mathscr{P}$ is the Weierstrass function. We consider the element
$$(p^{-n}\Omega_\infty)_{n=1}^\infty\in\varprojlim_n(p^{-n}\Omega_\infty F/\Omega_\infty F) = (\varprojlim_np^{-n}\mathscr{L}/\mathscr{L})\otimes\Q_p$$
which is in the Tate module of $\C/\mathscr{L}$ tensored with $\Q_p$. Let $V_pE_\omega = T_p E_\omega\otimes \Q_p$, $V_\p E_\omega = T_\p E_\omega\otimes\Q_p$, $V_{\overline{\p}}E_\omega = T_{\overline{\p}}E_\omega \otimes\Q_p$, and let $\xi = (\xi_n)_{n=1}^\infty$ be the image of $(p^{-n}\Omega_\infty)_{n=1}^\infty$ under the composition
$$\varprojlim_np^{-n}\Omega_\infty F/\Omega_\infty F\xrightarrow{\Phi_p}V_pE_\omega\longrightarrow V_\p E_\omega$$
where the second map is the projection corresponding to $T_pE_\omega = T_\p E_\omega\times T_{\overline{\p}} E_\omega$. 

The coordinates $x$ and $y$ on $E_\omega$ determine a formal group of $E$ over $K_\pa$, $\widehat{E}_\omega$. Let $V_p \widehat{E}_\omega = T_p\widehat{E}_\omega\otimes\Q_p$. Since $p$ splits in $F$ and $\p$ is the prime of $F$ determined by $\iota_p$, $T_p\widehat{E}_\omega = T_\p E_\omega$. Let $\xi$ now denote the corresponding elemet of $V_p\widehat{E}_\omega$. Since $V_p\widehat{E}$ is a rank one $\Q_p$-module, $\xi$ is a basis element. Let 
$$\zeta = (\zeta_{p^n})_{n = 1}^\infty = (\iota_p^{-1}(\exp(2\pi i/p^n)))_{n = 1}^\infty$$
so $\zeta$ is a basis element of $V_p\fG_m : = T_p\fG_m\otimes \Q_p$. Define
$$\varphi_p:V_p\widehat{E}_\omega\longrightarrow V_p\fG_m$$
by $\varphi_p(\xi) = \zeta$. It is a result of Tate (\cite{Ta2}) that the map
$$\Hom_{\OX_{\C_p}}(\widehat{E}_\omega,\fG_m)\longrightarrow \Hom_{\Z_p}(T_p\widehat{E}_\omega,T_p\fG_m)$$
is a bijection. We note that
$$\Hom_{\Q_p}(V_p\widehat{E}_\omega,V_p\fG_m) = \Hom_{\Z_p}(T_p\widehat{E}_\omega,T_p\fG_m)\otimes\Q_p$$
and let $\varphi\in \Hom_{\OX_{\C_p}}(\widehat{E},\fG_m)\otimes\Q_p$ be the element corresponding to $\varphi_p$. Define $\Omega_p$ by the rule
$$\omega = \Omega_p\varphi^*(dT/(1 + T)).$$

This defines a pair $(\Omega_\infty, \Omega_p)\in \C^\times\times\C_p^\times$. The definition depends on the choice of $E$, $\omega$, and $\gamma$, but is canonically defined as an element of $\C^\times\times\C_p^\times/\overline{\Q}^\times$. That is, if we make different choices for $E$, $\omega$, or $\gamma$, then $\Omega_\infty$ and $\Omega_p$ are both scaled by the same element of $\overline{\Q}^\times$ (see \cite{Ferr} for further explanation of the dependence).

\subsection{The CM Hida family}\label{CMHidaFamily}

For the remainder of Section \ref{ImagQuadChapter}, fix a nontrivial ray class character $\chi$ of conductor $\mathfrak{f}$ such that $(\mathfrak{f},p) = 1$, and let $f = \sum_{\ga\subset\OX_F}\chi(\ga)q^{N\ga}$ be the weight one modular form associated to $\chi$. Let $f_\alpha$ be a $p$-stabilization of $f$, so $\alpha$ is either $\chi(\p)$ of $\chi(\overline{\p})$. Recall that the character of $f$ is $\varepsilon:(\Z/N\Z)^\times\rightarrow\overline{\Q}$ determined by the rule $\varepsilon(\ell) = \chi(\ell\OX_F)$ for primes $\ell\nmid Np$. The goal of this section is to explicitly describe the rigid analytic functions $T_\ell$ for $\ell\nmid Np$ and $U_p$ on a neighborhood of the point corresponding to $f_\alpha$ on the eigencurve.

For $k\in\Z$, let $\nu_k\in\We(\Q_p)$ denote the character $t\mapsto t^{k-2}$. By Theorem \ref{BelDmit}), the eigencurve is \'{e}tale at the point corresponding to $f_\alpha$. Let $w = \nu_{1}\in\We(\Q_p)$ and let $W = W(w,1/p^r) = \Sp R$ be a neighborhood of $w$ such that the weight map $C_W^-\rightarrow W$ is \'{e}tale at all points in the connected component containing the point corresponding to $f_\alpha$. Let $x\in C_W^-(\C_p)$ be the point corresponding to $f_\alpha$ and let $V_{\C_p} = \Sp T_{\C_p}\subset C_{W,\C_p}^-$ be the connected component of $C_{W,\C_p}^-$ containing $x$. Then $V_{\C_p}\rightarrow W_{\C_p}$ is \'{e}tale, and we take $W$ to be as in Proposition \ref{BelProp}. Then the weight map on the level of rings $R_{\C_p}\rightarrow T_{\C_p}$ is an isomorphism, and we use this map to identify $T_{\C_p}$ with $R_{\C_p}$. 

Fix a choice of topological generator $\gamma$ of $1 + p\Z_p$, so
$$R = \left\{\sum a_n(t - (w(\gamma) - 1))^n \in\Q_p[[t - (w(\gamma) - 1)]]: |a_np^{rn}|\rightarrow 0\text{ as } n\rightarrow\infty\right\}.$$
Let $z = t - (w(\gamma) - 1)$. Then $R$ is the set of all $F(z) \in\Q_p[[z]]$ that converge on the closed around $0$ disk of radius $1/p^r$ in $\C_p$. By the Weierstrass preparation theorem, any $F(z)\in R$ is determined by its values
$$ev_{\nu_k}(F(z)) = F(\nu_k(\gamma) - w(\gamma)) = F(\gamma^{k-2} - \gamma^{-1})$$
at the integers $k\in\Z$ such that $\nu_k\in W$.  For an integer $k$, $\nu_{k}$ is in  $W = W(w,1/p^r)$ if and only if $k\equiv1 \bmod p^{r-1}(p-1)$.

Since $V$ is \'{e}tale over weight space, the Hecke operators $T_\ell$ for $\ell\nmid Np$ $U_p$, and $[a]$ for $a\in(\Z/N\Z)^\times$ as rigid analytic functions in $R_{\C_p}$ are deterined by the following two properties: 
\begin{enumerate}[nolistsep]
\item At the weight $w$,
$$ev_w(T_\ell) = a_\ell(f_\alpha) = \begin{cases} \chi(\mathfrak{q}) + \chi(\overline{\mathfrak{q}}) &\text{ if } \ell\OX_F = \mathfrak{q}\overline{\mathfrak{q}}\\ 0 & \text{ if } \ell \text{ is inert in }F\end{cases}$$
$ev_w(U_p) = \alpha$, and $ev_w([a]) = \varepsilon(a)$ for all $a\in (\Z/N\Z)^\times$.

\item For all $k\in\Z_{\geq 2}$ such that $\nu_k\in W$, $ev_{\nu_k}(T_\ell), ev_{\nu_k}(U_p)$ are the $T_\ell$ and $U_p$ Hecke eigenvalues of an eigenform $g$ of weight $k$, level $\Gamma_0$, and character $\varepsilon$ which is new at level $N$.
\end{enumerate}

The second condition implies that the functions $[a]\in R_{\C_p}$ are the constant function $[a] = \varepsilon(a)$. We exhibit explicit elements of $R_{\C_p}$ with the above two properties as $T_\ell$ for $\ell\nmid Np$ and $U_p$. 

In the interest of clarity of composition and space, we assume for the rest of Section \ref{ImagQuadChapter} that $\alpha = \chi(\overline{\p})$. The case $\alpha = \chi(\p)$ is similar (see \cite{Ferr} for more details).

To begin we define an algebraic Hecke character of $F$. Since $p\geq 3$, the only root of unity (and so the only unit of $F$) congruent to $1\bmod \p$ in $F$ is $1$. Therefore we may identify the group $P_1(\p)$ with a subgroup of $F^\times$:
$$P_1(\p) = \{\alpha\in F^\times : ((\alpha), \p) = 1, \alpha\equiv 1\bmod \p\}\subset F^\times.$$
Define $\lambda_0$ as
$$\lambda_0:P_1(\p)\longrightarrow F^\times\subset\overline{\Q}^\times$$
$$\lambda_0(\alpha) = \alpha = \sigma_1(\alpha).$$
Since $\overline{\Q}^\times$ is divisible, we may extend $\lambda_0$ to $I(\p)$ to define an algebraic Hecke character $\lambda$ of infinity type $(1,0)$ and modulus $\p$. The choice of extension of $\lambda_0$ is determined up to multiplication be a character of $I(\p)/P_1(\p)$. We impose a condition on the extension $\lambda$ we choose. Recall that $\C_p^\times$ may be written as $\C_p^\times =  p^\Q\times W\times U$, where $W$ is the group of roots of unity of order prime to $p$ and  $U = \{u\in\C_p^\times : |1 - u|<1\}$. By construction, after composing with $\iota_p$ the image of $\lambda_0$ is contained in $U$. Since $U$ is a divisible group, we may choose our extension $\lambda$ so that the image of $\lambda$ after composing with $\iota_p$ is also contained in $U$, which we do.  Since the only torsion elements in $U$ are the $p$-power roots of unity, any two extensions $\lambda$ and $\lambda'$ of $\lambda_0$ that have image in $U$ differ by a character of $I(\p)/P_1(\p)[p^\infty]$ where the $[p^\infty]$ denotes the maximal quotient of $I(\p)/P_1(\p)$ with $p$-power order. 

Let $p^n = |I(\p)/P_1(\p)[p^\infty]|$. If $p^r \leq p^n$, then we shrink $W$ so that $W = W(w, \frac{1}{p^{n + 1}})$. We may do this without changing anything we have assumed previously, and the reason for doing this will become clear momentarily. 

Let $M = |I(\p)/P_1(\p)|$ and note that $|M|_p = 1/p^n$. For each prime $\q$ of $F$ such that $\q\not=\p$ define the power series
$$G_\q(z) = \exp_p(z\log_p(\lambda(\q))) = \sum_{n = 0}^\infty \frac{z^n\log_p(\lambda(\q))^n}{n!}$$
as an element of $\C_p[[z]]$. The power series $G_\q(z)$ converges if
$$|z|<\frac{1}{p^{1/(p-1)}|\log_p(\lambda(\q))|}.$$
Since $M = |I(\p)/P_1(\p)|$, $\q^M = (q)$ for some $q\in\OX_F$ such that $q\equiv 1\bmod \p$. Hence by definition of $\lambda_0$
$$\lambda(\q)^M = \lambda((q))\equiv 1\bmod \p$$
so $|1 - \lambda(\q)^M| < p^{-1/(p-1)}$. Then by properties of the $p$-adic logarithm,
$$\frac{1}{p^{1/(p-1)}}>|1 - \lambda(\q)^M| = |\log_p(\lambda(\q)^M)| = |M||\log_p(\lambda(\q))| = \frac{|\log_p(\lambda(\q))|}{p^n}$$
so
$$\frac{1}{p^n}<\frac{1}{p^{1/(p-1)}|\log_p(\lambda(\q))|}.$$
Therefore $G_\q(z)$ converges for $|z|\leq\frac{1}{p^n}$, which is independent of $\q$.

Recall that $\log_\gamma(z):=\frac{\log_p(z)}{\log_p(\gamma)}$, and define
$$F_\q(z) = G_\q\circ\log_\gamma(1 + \gamma z).$$
By construction, if $|z|\leq \frac{1}{p^{n+1}}$ then $F_\q(z)$ converges. This implies that $F_\q(z) \in R_{\C_p}$. The function $F_\q(z)$ is the unique element of $R_{\C_p}$ with the property that for all $k\in\Z$ such that $\nu_k\in W$, $ev_{\nu_k}(F_\q(z)) = (\lambda(\q))^{k-1}$. Furthermore, since $k\in\Z$ is such that $\nu_k\in W$ if and only if $k\equiv 1\mod p^{r-1}(p-1)$ and $r>n$, $F_\q(z)$ does not depend on the choice of extension $\lambda$ of $\lambda_0$ since $p^n$ divides $k-1$ so the exponent $k-1$ will kill any character of $I(\p)/P_1(\p)[p^\infty]$. 

Now let $\ga\subset\OX_F$ be a nontrivial ideal of $\OX_F$ such that $(\ga,\p) = 1$, and define
$$F_\ga(z) = \begin{cases} \ds\prod_{\q} F_\q(z)^{val_\q(\ga)} &\text{if } (\ga,\p) = 1\\
						0 &\text{else.}\end{cases}$$
Further, define $A_1(z) = 1$ and for $n\geq 2$ define
$$A_n(z) = \sum_{\substack{\ga\subset\OX_F\\ N_{F/\Q}\ga = n}}\chi(\ga)F_\ga(z).$$
Define the formal $q$-expansion
$$\mathcal{F} = \sum_{n = 1}^\infty A_n(z)q^n\in R_{\C_p}[[q]].$$
This formal $q$-expansion is the CM Hida family specializing to $f_\alpha$ in weight one.
\begin{pr} For all $k\in\Z_{\geq1}$, $\nu_k\in W$
$$\mathcal{F}_k : = \sum_{n = 1}^\infty ev_{\nu_k}(A_n(z))q^n = \sum_{\ga\subset\OX_F}\chi\lambda^{k-1}(\ga)q^{N\ga}$$
is the $q$-expansion of a weight-$k$ cusp form of level $\Gamma_0$ and character $\varepsilon$ that is new at level $N$.
\end{pr}

\begin{proof}
By definition of $A_n(z)$ we have that
$$\sum_{n = 1}^\infty ev_{\nu_k}(A_n(z))q^n = \sum_{\ga\subset\OX_F}\chi\lambda^{k-1}(\ga)q^{N\ga}.$$
Shimura (\cite{Sh1}) showed that 
$$\sum_{\ga\subset\OX_F}\chi\lambda^{k-1}(\ga)q^{N\ga}$$
is the $q$-expansion of a weight-$k$ cusp form of level $\Gamma_0$ which is new at level $N$ and has character defined by 
$$\ell \longmapsto \frac{\chi((\ell))\lambda^{k-1}((\ell))}{\ell^{k-1}} = \chi((\ell))\left(\frac{\lambda((\ell))}{\ell}\right)^{k-1}$$
for $\ell\in(\Z/N\Z)^\times$ a prime not equal to $p$. A simple calculation shows that this is the character $\varepsilon$.
\end{proof}

By the proposition, the functions $A_\ell(z)\in R_{\C_p}$ for $\ell\nmid Np$ and $A_p(z)\in R_{\C_p}$ satisfy the two properties that uniquely determine $T_\ell, U_p\in R_{\C_p}$. Hence $T_\ell = A_\ell$ for $\ell\nmid Np$ and $U_p = A_p$.

\subsection{Two-variable $p$-adic $L$-function of the CM family}\label{TwoVarCM}

Keeping the notation of the previous section, let $\Phi$ be a generator for the rank one $T_{\C_p}$-module 
$$\Symb_{\Gamma_0}^-(\D(R))^o\otimes_{\T_{W}^-} T_{\C_p}\subset \Symb_{\Gamma_0}^-(\D(R))^o$$
and let
$$\Lambda(\Phi,\cdot,\cdot): W\times\We\longrightarrow \C_p$$
be the two-variable $p$-adic $L$-function associated to $\Phi$ as in Section \ref{TwoVariable}. In order to prove Conjecture \ref{pAdicStarkAt0} we restrict $\Lambda(\Phi,\cdot,\cdot)$ to a particular subset of $W\times\We$. Let
$$U = \{t\in\Z_p : t\equiv 1\bmod p^{r-1}\}$$
where $W = W(w,\frac{1}{p^r})$ and $r$ was chosen in the previous section. Let $\eta$ be a $p$-power order character and let $\psi = \eta\omega$. Let $p^m$ be the conductor of $\psi$. Define the two-variable restriction of $\Lambda(\Phi,\cdot,\cdot)$:
$$L_p(\chi\eta\omega,\alpha,\cdot,\cdot):U\times\Z_p\longrightarrow \C_p$$
$$L_p(\chi\eta\omega,\alpha,t,s) = \Lambda(\Phi,\omega^{-1}\langle\cdot\rangle^{t-2},(\eta\omega)^{-1}\langle\cdot\rangle^{s-1}).$$
For all $k\in\Z_{\geq 2}, k\equiv 1\bmod p^{r-1}$, let $(\Omega_{\infty,k},\Omega_{p,k})$ be the periods for $\nu_{k}\in W$ that appear in the interpolation formula for $\Lambda(\Phi,\cdot, \cdot)$. Then $L_p(\chi\eta\omega,\alpha, t, s)$ is determined by the following interpolation property: for all $k\in\Z_{\geq 2}$, $k\equiv 1\bmod p^{r-1}$, and $j\in \Z$, $1\leq j\leq k-1$, $j\equiv 1\bmod 2(p-1)$
$$\frac{L_p(\chi\eta\omega,\alpha,k,j)}{\Omega_{p,\nu_k}} = E_p(\alpha,\eta\omega, k, j)\frac{L(\chi\lambda^{k-1}\eta\omega,j)}{\Omega_{\infty,\nu_k}}$$
where
\[\arraycolsep=1.4pt\def\arraystretch{2.2}
\begin{array}{ll} \ds E_p(\alpha,\eta\omega,k,j) & = \ds \ds\frac{1}{\chi\lambda^{k-1}(\overline{\p})^m}\left(1 - \frac{(\eta\omega)^{-1}(p)p^{j-1}}{\chi\lambda^{k-1}(\overline{\p})}\right)\times\\ &\ds \times\frac{p^{m(j-1)} (j-1)!\tau((\eta\omega)^{-1})}{(2\pi )^{j-1}}\end{array}\]
and $L(\chi\lambda^{k-1}\eta\omega,s)$ is the complex Hecke $L$-function associated to $\chi\lambda^{k-1}\eta\omega$.

\subsection{Two-variable specialization of $L_{p,Katz}$}

In this section we define a two-variable specialization of Katz's $p$-adic $L$-function that we compare to the two-variable $p$-adic $L$-function defined in the previous section.

Observe that the complex $L$-value appearing the interpolation formula in the previous section is
$$L(\chi\lambda^{k-1}\eta\omega,j) = L(\chi\lambda^{k-1}\eta\omega N^{-j},0).$$
By our choice of $\lambda$, the algebraic Hecke character $\chi\lambda^{k-1}\eta\omega N^{-j}$ has infinity type $(k-1-j,-j)$, which is not in the range of interpolation of Katz's $p$-adic $L$-function.

From here on, let $c$ denote complex conjugation, so $c$ is an automorphism of $\C$. Via our embedding $\iota_\infty$, $c$ acts on ideals of $F$, and there is the relation of complex $L$-functions
$$L(\chi\lambda^{k-1}\eta\omega N^{-j},s) = L(\chi\lambda^{k-1}\eta\omega N^{-j}\circ c,s)$$
that changes the infinity type. Therefore, $\chi\lambda^{k-1}\eta\omega N^{-j}\circ c$ has infinity type $(-j, k-1-j)$, which is in the range of interpolation of Katz's $p$-adic $L$-function. 

Let $\kappa_1 = \lambda\circ c$ viewed as an algebraic Hecke character. By our choice of $\lambda$, $\kappa_1$ has infinity type $(0,1)$ and conductor $\overline{\p}$. Further when we view $\kappa_1$ as a $p$-adic Hecke character, since $\lambda$ takes values in $U = \{u\in\C_p^\times : |1-u|<1\}\subset\C_p^\times$ we may consider the $p$-adic Hecke character $\kappa_1^{s_1}$ for any $p$-adic number $s_1\in\Z_p$. 

Let $\kappa_2$ be the algebraic Hecke character $\kappa_2 = \omega^{-1}N$ where $N$ is the norm character
$$N:\A_F\longrightarrow \overline{\Q}^\times$$
$$N((x_v)_v) = \prod_{v-\text{finite}} |x_v|^{-1}.$$
Viewing $\kappa_2$ as a $p$-adic Hecke character, $\kappa_2$ has image $1 + p\Z_p$ in $\C_p^\times$. It therefore makes sense to consider $\kappa_2^{s_2}$ as a $p$-adic Hecke character for any $s_2\in\Z_p$. Let $\widetilde{\chi} = \chi\circ c$ and note that $\widetilde{\chi}$ has conductor $\overline{\mathfrak{f}}$. Let $L_{p,Katz}$ be Katz's $p$-adic $L$-function with respect to the ideal $\mathfrak{m}$ where as in the notation of Section \ref{DefinitionsOfFunctions}, $\mathfrak{m}$ is the conductor of $M/F$. The ideal $\mathfrak{m}$ is divisible by all the primes that divide $\mathfrak{f}$ and $\overline{\mathfrak{f}}$. Let $(\Omega_\infty,\Omega_p)$ be the period pair used to define $L_{p,Katz}$. 

Define
$$L_{p,Katz}(\chi\eta,\alpha,\cdot,\cdot):U\times\Z_p\longrightarrow \C_p$$
$$L_{p,Katz}(\chi\eta,\alpha,s_1,s_2) : = L_{p,Katz}(\widetilde{\chi}\eta \kappa_1^{s_1 - 1}\kappa_2^{-s_2}).$$

\begin{pr} $L_{p,Katz}(\chi\eta,\alpha,s_1,s_2)$ is determined by the following interpolation property: for all $k\in\Z_{\geq 2}, k\equiv 1\bmod p^{r-1}$, $j\in \Z$, $1\leq j\leq k-1$, $j\equiv 1\bmod p-1$,
$$\frac{L_{p,Katz}(\chi\eta,\alpha,k,j)}{\Omega_p^{k-1}} = E_{p}(\alpha,\eta\omega,k,j) \frac{-(2\pi  )^{k-2}}{\sqrt{d_F}^{k-1-j}}\frac{L(\chi\lambda^{k-1}\psi\omega^{j-1},j)}{\Omega_\infty^{k-1}}$$
where
$E_p(\alpha,\eta\omega,k,j)$ is defined as in the previous section. \end{pr}

\begin{proof} That $L_{p,Katz}(\chi\eta,\alpha,s_1,s_2)$ is determined by the interpolation property follows from the continuity of $L_{p,Katz}(\chi\eta,\alpha,s_1,s_2)$ and that the set of $k$'s and $j$'s is dense in $U\times\Z_p$. Let $k\in\Z_{\geq 2}$, $k\equiv 1\bmod p^{r-1}$ and $j\in\Z$, $1\leq j\leq k-1$, $j\equiv 1\bmod p-1$. By our definitions
\begin{equation} \widetilde{\chi}\eta\kappa_1^{k-1}\kappa_2^{-j} =  \chi\eta\omega\lambda^{k-1}N^{-j}\circ c \end{equation}
so $\widetilde{\chi}\eta\kappa_1^{k-1}\kappa_2^{-j}$ has infinity type $(-j, k-1-j)$ which is in the range of interpolation for $L_{p,Katz}$. 

By the interpolation formula for $L_{p,Katz}$, 
\[\arraycolsep=1.4pt\def\arraystretch{2.2}
\begin{array}{ll} \ds\frac{L_{p,Katz}(\chi\eta,\alpha,k,j)}{\Omega_p^{k-1}} 
	& \ds = \frac{(j-1)!(2\pi)^{k-1-j}}{\sqrt{d_F}^{k-1-j}}W_p(\chi\eta\omega\lambda^{k-1}N^{-j}\circ c)\times\\
	& \ds \times \left(1 - \frac{(\chi\eta\omega)^{-1}\lambda^{1-k}N^j(\overline{\p})}{p}\right)(1 - \chi\eta\omega\lambda^{k-1}N^{-j}(\p))\times\\
	& \ds \times\frac{L(\chi\eta\omega\lambda^{k-1}N^{-j},0)}{\Omega_\infty^{k-1}}
\end{array}\]
Since $\lambda$ has modulus $\p$, $1 - \chi\eta\omega\lambda^{k-1}N^{-j}(\p) = 1$. We also have that $(\eta\omega)^{-1}(\overline{\p}) = (\eta\omega)^{-1}(p)$, $N^j(\overline{\p}) = p^j$, and a calculation shows that, 
$$W_p(\chi\eta\omega\lambda^{k-1}N^{-j}\circ c) = \frac{-p^{m(j-1)}\tau((\eta\omega)^{-1})}{\chi\lambda^{k-1}(\overline{\p})^m}.$$
Therefore the formula becomes
\[\arraycolsep=1.4pt\def\arraystretch{2.2}
\begin{array}{ll} \ds\frac{L_{p,Katz}(\chi\eta,\alpha,k,j)}{\Omega_p^{k-1}} 
	& \ds = \frac{(j-1)! (2\pi)^{k-1-j}}{\sqrt{d_F}^{k-1-j}}\frac{-p^{m(j-1)}\tau((\eta\omega)^{-1})}{\chi\lambda^{k-1}(\p)^m}\times\\
	& \ds \times \left(1 - \frac{(\eta\omega)^{-1}(p)p^{j-1}}{\chi\lambda^{k-1}(\overline{\p})^m}\right)\frac{L(\chi\lambda^{k-1}\eta\omega,j)}{\Omega_\infty^{k-1}}\\
	& \ds= E_{p}(\alpha,\eta\omega,k,j) \frac{-(2\pi )^{k-2}}{\sqrt{d_F}^{k-1-j}}\frac{L(\chi\lambda^{k-1}\psi\omega^{j-1},j)}{\Omega_\infty^{k-1}}.\end{array}\]\end{proof}

\subsection{Choice of periods and comparison}\label{ChoiceOfPeriodsAndComparison}

Let $S_{\C_p}$ be the fraction field of $T_{\C_p} = R_{\C_p}$. 

\begin{pr} There exists $\Psi\in \Symb_{\Gamma_0}^-(\D(R))\otimes_{\T_W^\pm}S_{\C_p}$ such that the $p$-adic $L$-function
$$L_p(\chi\eta\omega,\alpha,t,s) : = \Lambda(\Psi, \omega^{-1}\langle\cdot\rangle^{t-2},(\eta\omega)^{-1}\langle\cdot\rangle^{s-1})$$
is calculated with the $p$-adic and complex periods
$$(\Omega_{p,k},\Omega_{\infty,k}) = \left(\Omega_p^{k-1},\Omega_\infty^{k-1}\left(\frac{\sqrt{d_F}}{2\pi}\right)^{k-2}\right)$$
where $(\Omega_p,\Omega_\infty)$ is the period pair used to define Katz's $p$-adic $L$-function. We note that the domain of $L_p(\chi\eta\omega,\alpha,t,s)$ is as in the previous section.\end{pr}
\begin{proof} Let $L_p(\chi\eta\omega,\alpha,t,s) = \Lambda(\Phi,\omega^{-1}\langle\cdot\rangle^{t-2},(\eta\omega)^{-1}\langle\cdot\rangle^{s-1})$ be as in Section \ref{TwoVarCM}. We determine a meromorphic function $P(t)$ on $U$ such that $P(t)L_p(\chi\eta\omega,\alpha,t,s)$ has interpolation formula with the periods
$$\left(\Omega_p^{k-1},\Omega_\infty^{k-1}\left(\frac{\sqrt{d_F}}{2\pi}\right)^{k-2}\right).$$

Let 
$$P:U\times \Z_p\longrightarrow \C_p\cup\{\infty\}$$
be the $p$-adic meromorphic function defined by the ratio
$$P(t,s) = \frac{L_{p,Katz}(\chi\eta,\alpha,t,s)}{L_p(\chi\eta\omega,\alpha,t,s)}.$$
Then $P(t,s)$ has the interpolation property:
$$\frac{P(k,j)\Omega_{p,k}}{\Omega_p^{k-1}} = \frac{\Omega_{\infty,k}}{\Omega_\infty^{k-1}}\frac{-(2\pi)^{k-2}}{\sqrt{d_F}^{k-1-j}}$$
for $k$'s and $j$'s as in the previous section. 

When defining the periods for $L_p(\chi\eta\omega,\alpha,t,s)$, we choose $\Phi$ which we've done, and we choose the $\Omega_{\infty,k}$. These choices determine the $\Omega_{p,k}$. The condition on the choice of $\Omega_{\infty,k}$ is that the complex values in the interpolation formula $L_p(\mathcal{F}_k,\cdot,\cdot)$ are algebraic. These values are for all odd finite order characters $\psi\in \We(\C_p)$, $k\in\Z_{\geq 2}$, $1 \leq j\leq k-1$,
$$C_{alg}(\alpha,k,j)\frac{L(\chi\lambda^{k-1}\psi\omega^{j-1},j)}{(2\pi )^{j-1}\Omega_{\infty,k}}$$
where
$$C_{alg}(\alpha,k,j) = \ds\frac{p^{m(j-1)}(j-1)!\tau(\psi^{-1}\omega^{1-j})}{\chi\lambda^{k-1}(\overline{\p})^m}\left(1 - \frac{\psi^{-1}\omega^{1-j}(p)}{\chi\lambda^{k-1}(\overline{\p})p^{1-j}}\right)\frac{1}{i^{j-1}}$$
and $m$ is the power of $p$ in the conductor of $\psi$.

It is clear then that we may define
$$\Omega_{\infty,k} = \Omega_\infty^{k-1}\left(\frac{\sqrt{d_F}}{2\pi }\right)^{k-2}$$
since by the interpolation property of Katz's $p$-adic $L$-function, the values
$$C_{alg}(\alpha,k,j)\frac{(2\pi )^{k-1-j} L(\chi\lambda^{k-1}\psi\omega^{j-1},j)}{\sqrt{d_F}^{k-2}\Omega_\infty^{k-1}}$$
are algebraic.

If we consider $P(t,s)$ with this choice of complex periods, then $P(t,s)$ satisfies the interpolation formula for $k\in\Z_{\geq 2}$ $k\equiv 1\bmod p^{r-1}$, $j\in\Z$, $1\leq j\leq k-1$, $j\equiv 1\bmod 2(p-1)$
$$\frac{P(k,j)\Omega_{p,k}}{\Omega_p^{k-1}} = -\sqrt{d_F}^{j-1}.$$
Now we separate variables for the function $P(t,s)$. Since $p$ splits in $F$, $\sqrt{d_F}\in\Q_p=F_\p$. Define the analytic function $Q(s)$ as $Q(s) = -\langle \sqrt{d_F}\rangle^{s-1}$, and let $P(t) = P(t,s)/Q(s)$. The function $P(t)$ is a $p$-adic meromorphic function on $U$ satisfying the relation that for all $k\in\Z_{\geq 2}$, $k\equiv 1\bmod p^{r-1}$, $P(k)\Omega_{p,k} = \Omega_p^{k-1}$. Since $P(t)$ is a $p$-adic meromorphic function on $U$, there exists an element $\widetilde{P}\in S_{\C_p}$ such that for all $t\in U$, $\widetilde{P}(\gamma^{t-2} - \gamma^{-1}) = P(t)$.
 
If we define $\Psi = \widetilde{P}\Phi^-$ and redefine the function
$$L_p(\chi\eta\omega,\alpha,t,s) = \Lambda(\Psi, \omega^{-1}\langle\cdot\rangle^{t-2},(\eta\omega)^{-1}\langle\cdot\rangle^{s-1})$$
then $L_p(\chi\eta\omega,\alpha,t,s)$ satisfies the interpolation property that for all $k,j$ as above,
$$\frac{L_p(\chi\eta\omega,\alpha,t,s)}{\Omega_p^{k-1}} = E_p(\alpha,\eta\omega,k,j)\frac{(2\pi)^{k-2}L(\chi\lambda^{k-1}\psi\omega^{j-1},j)}{\sqrt{d_F}^{k-2}\Omega_\infty^{k-1}}.$$
\end{proof}

\begin{rem} If $P(t)$ in the proof of the previous proposition does not have any zeros or poles, then $\Psi$ is a generator for the free rank one $T_{\C_p}$-module $\Symb_{\Gamma_0}(\D(R))^o\otimes_{\mathbb{T}_W^-}T_{\C_p}$ and so $\Psi$ would be a valid choice to define the $p$-adic $L$-function as in Section \ref{TwoVariable}.\end{rem}

We record the precise comparison of the $p$-adic $L$-function defined in the previous two sections that appeared in the proof of the previous proposition.

\begin{cor}\label{PeriodsCorollary} Let $L_{p,Katz}(\chi\eta,\alpha,t,s)$ and $L_p(\chi\eta\omega,\alpha,t,s)$ be defined as in the previous two sections, so 
$$L_p(\chi\eta\omega,\alpha,t,s) = \Lambda(\Phi,\omega^{-1}\langle\cdot\rangle^{t-2},(\eta\omega)^{-1}\langle\cdot\rangle^{s-1})$$
where $\Phi$ is a generator of $\Symb_{\Gamma_0}^-(\D(R))^o\otimes_{\T_{W}^\pm}T_{\C_p}$ as a $T_{\C_p}$-module. Then 
$$L_{p,Katz}(\chi\eta,\alpha,t,s) = P(\eta,t,s) L_p(\chi\eta\omega,\alpha,t,s)$$
where $P(\eta,t,s)$ is a $p$-adic meromorphic function determined by the interpolation property that for all $k\in\Z_{\geq 2}$, $k\equiv 1\bmod p^{r-1}$, $j\in\Z$, $1\leq j\leq k-1$, $j\equiv 1\bmod 2(p-1)$,
$$\frac{P(\eta,k,j)\Omega_{p,k}}{\Omega_p^{k-1}} = \frac{\Omega_{\infty,k}}{\Omega_\infty^{k-1}}\frac{-(2\pi)^{k-2}}{\sqrt{d_F}^{k-1-j}}.$$
\end{cor}

\begin{rem} 
We remark that $P(\eta,t,s)$ a priori depends on $\eta$ and $\alpha$, but as is clear from the interpolation formula does not actually depend on $\eta$ or $\alpha$. The reason for putting $\eta$ in the notation will become clear in the next section.\end{rem}

\subsection{Proof of the conjecture in this case}\label{finishProofKatz}

In this section we prove Conjecture \ref{pAdicStarkAt0} for $\chi$. We adopt the notation of Section \ref{DefinitionChapter}. For each $r\geq 1$ let $u_r = u_{M_r}$ be the Stark unit for $M_r/F$ from Definition \ref{ComplexStarkUnitDef}. For $\varphi\in \We(\C_p)$ a character of order $p^r$, the unit $u_{\chi\varphi,\alpha}^*$ is obtained from $u_r$ by first mapping $u_r$ to the $(\rho\varphi)^*$-isotypic component of $k_r\otimes U_r$ and then projecting to the subspace where $\delta_p$ acts with eigenvalue $\alpha$. Let $\pi_{\rho\varphi}^*$ be the map
$$\pi_{\rho\varphi}^*:U_r\longrightarrow (k_r\otimes U_r)^{(\rho\varphi)^*}$$
$$\pi_{\rho\varphi}^*(u) = \sum_{\sigma\in \Delta_n} Tr((\rho\varphi)^*(\sigma))\otimes \sigma(u).$$
The idea to project to the subspace where $\delta_p$ acts with eigenvalue $\alpha$ is of Greenberg and Vatsal (\cite{GV}) and we adopt their notation. Let $|\cdot|_\alpha$ denote the map
$$|\cdot|_\alpha:(k\otimes U_r)^{(\rho\varphi)^*}\longrightarrow (k_r\otimes U_r)^{(\rho\varphi)^*, \delta_p = \alpha}$$
$$|u|_\alpha = \frac{1}{|\Delta_p|}\sum_{i = 0}^{|\Delta_p| - 1}\alpha^{-i}\delta_p^i(u).$$
Then $|\pi_{\rho\varphi}^*(u_r)|_{\alpha}\in (k_r\otimes U_r)^{(\rho\varphi)^*, \delta_p = \alpha}$ and so the following theorem implies Conjecture \ref{pAdicStarkAt0}.

\begin{thm}\label{finishThmKatz} Let $\eta,\psi\in\We(\C_p)$ be of orders $p^m$ and $p^n$ respectively. Then
$$L_p(\chi,\alpha,\eta\omega,\psi\omega,0) = \frac{\frac{\tau(\eta^{-1})}{p^{m+1}}\left(1 - \frac{\eta^{-1}(p)}{\alpha p}\right)(1 - \beta\eta(p))}{\frac{\tau(\psi^{-1})}{p^{n+1}}\left(1 - \frac{\psi^{-1}(p)}{\alpha p}\right)(1 - \beta\psi(p))}\frac{\log_p|\pi_{\rho\eta}^*(u_m)|_{\alpha}}{\log_p|\pi_{\rho\psi}^*(u_n)|_{\alpha}}.$$
\end{thm}

\begin{proof} To begin, we simplify the expression $|\pi_{\rho\eta}^*(u_m)|_{\alpha}$. Since $(\rho\eta)^* = \Ind_{H_m}^{\Delta_m}(\chi\eta)^{-1}$, for all $\sigma\in \Delta_m - H_m$, $Tr((\rho\eta)^*(\sigma)) = 0$. Since $c\in \Delta_m - H_m$, for all $\sigma\in H_m$, $Tr((\rho\eta)^*(\sigma)) = \chi\eta(\sigma) + \chi\eta(c\sigma c)$. Let $\chi_c$ denote the character $\chi_c(\sigma) = \chi(c\sigma c)$ and note that since $\Q_n$ is totally real, $\eta(c\sigma c) = \eta(\sigma)$ for all $\sigma$. Therefore,
$$\pi_{\rho\eta}^*(u_m) = \sum_{\sigma\in H_m} \chi\eta(\sigma)\otimes\sigma(u_m) + \chi_c\eta(\sigma)\otimes\sigma(u_m).$$
Since $\alpha = \chi(\overline{\p})$, we have that
$$\lvert \sum_{\sigma\in H_m}\chi\eta(\sigma)\otimes\sigma(u_m)\rvert_{\alpha} = 0 \text{ and } \lvert \sum_{\sigma\in H_m} \chi_c\eta(\sigma)\otimes\sigma(u_m)\rvert_{\alpha} = \sum_{\sigma\in H_m} \chi_c\eta(\sigma)\otimes\sigma(u_m).$$
Therefore
$$|\pi_{\rho\eta}^*(u_m)|_{\alpha} = \ds\sum_{\sigma\in H_m} \chi_c\eta(\sigma)\otimes\sigma(u_m).$$
A similar formula holds for $|\pi_{\rho\psi}^*(u_n)|_{\alpha}$.

Let $L_p(\chi\eta\omega,\alpha,t,s)$ and $L_p(\chi\psi\omega,\alpha,t,s)$ be as defined in Section \ref{TwoVarCM}. By construction 
$$L_p(\chi,\alpha,\eta\omega,\psi\omega,s) = \frac{L_p(\chi\eta\omega,\alpha,1,s)}{L_p(\chi\psi\omega,\alpha,1,s)}.$$
Then by Corollary \ref{PeriodsCorollary},
\[\niceArray
\begin{array}{ll} L_p(\chi,\alpha,\eta\omega,\psi\omega,s) & = \ds\frac{P(\eta,1,s)L_{p,Katz}(\chi\eta\omega,\alpha,1,s)}{P(\psi,1,s)L_{p,Katz}(\chi\psi\omega,\alpha,1,s)}\\
				& \ds=\frac{L_{p,Katz}(\chi\eta\omega,\alpha,1,s)}{L_{p,Katz}(\chi\psi\omega,\alpha,1,s)}.\end{array}\]
Plugging in $0$, we get
$$L_p(\chi,\alpha,\eta\omega,\psi\omega,0)  = \ds\frac{L_{p,Katz}(\chi\eta\circ c)}{L_{p,Katz}(\chi\psi\circ c)}.$$
We now use Theorem \ref{KatzpAdicKron}. By the above simplifications of $|\pi_{\rho\eta}^*(u_m)|_{\alpha}$ and $|\pi_{\rho\psi}^*(u_n)|_{\alpha}$,
$$\frac{L_{p,Katz}(\chi\eta\circ c)}{L_{p,Katz}(\chi\psi\circ c)} = \frac{\ds\frac{\tau(\eta^{-1})}{p^{m+1}}\left(1 - \frac{(\chi\eta)^{-1}(\overline{\p})}{p}\right)(1 - \chi\eta(\p))\log_p|\pi_{\rho\eta}^*(u_m)|_{\alpha}}{\ds\frac{\tau(\psi^{-1})}{p^{n+1}}\left(1 - \frac{(\chi\psi)^{-1}(\overline{\p})}{p}\right)(1 - \chi\psi(\p))\log_p|\pi_{\rho\psi}^*(u_n)|_{\alpha}}.$$
To finish, we just note that since $\alpha = \chi(\overline{\p})$, $\beta = \chi(\p)$, so $(\chi\eta)^{-1}(\overline{\p}) = \eta^{-1}(p)/\alpha$ and $(\chi\psi)^{-1}(\p) = \psi^{-1}(p)/\alpha$, as well as $\chi\eta(\p) = \beta\eta(p)$ and $\chi\psi(\p) = \beta\psi(p)$.\end{proof}

\section{Numerical Evidence}\label{NumEvChap}

The programming for the examples consisted of three basic parts: computing the minimal polynomial of the Stark units, viewing the Stark units $p$-adically to take their $p$-adic logarithm, and computing the $p$-adic $L$-values. The code used for the examples can be found at \url{https://github.com/Joe-Ferrara/p-adicStarkExamples}. We briefly explain the basis of the code.

In the case where $F$ is real quadratic, the minimal polynomial of the Stark units was computed in SAGE combining the strategies of Stark in \cite{St1} and Dummit, Sands, and Tangedal in \cite{DST}. In the cases where $F$ is imaginary quadratic, the minimal polynomials of the Stark units were computed in pari/gp using the formulas from Section \ref{ComplexStarkUnits}. To view the Stark units $p$-adically and take the $p$-adic logarithm we wrote a class in SAGE to represent the extension of $\Q_p$ the Stark units are in and to take their $p$-adic logarithm. To compute the $p$-adic $L$-values, we used code written in SAGE by Rob Harron and Rob Pollack to compute overconvergent modular symbols (their code is based off the algorithms described in \cite{PS2}). We computed the weight negative one overconvergent modular symbol associated to $f_\alpha$ to get the $p$-adic $L$-values as described in Remark \ref{oneVarDef}.

An important reason for these examples is that we expected the units appearing in Conjecture \ref{pAdicStarkAt0} to be related to the Stark units in definition \ref{StarkUnitDef} in the way that they are related in Section \ref{finishProofKatz}, when $F$ in imaginary quadratic and $p$ is split in $F$ (see \cite{Ferr} for what we expected). As the examples show this may be the case when $F$ is real quadratic. When $F$ is imaginary quadratic and $p$ is inert in $F$, we can verify the conjecture, but it is not clear how or if the units in Conjecture \ref{pAdicStarkAt0} are related to the Stark units in \ref{ComplexStarkUnitDef}. In the case when $F$ is imaginary quadratic and $p$ is inert in $F$, the expected formulas conjectured in \cite{Ferr} are not correct.

We adopt all the notation of Section \ref{DefinitionChapter}. All three examples are of the following form which we describe before specifying the exact examples. 

Let $\psi\in\We(\C_p)$ be the character $\psi:(\Z/p^2\Z)^\times\rightarrow\overline{\Q}^\times$ that sends the generator of $(\Z/p^2\Z)^\times$ with minimal positive integer coset representative, to $\zeta_p$. For $\alpha = \pm 1$ in the first example and $\alpha = -1$ in the second two examples, we verify the conjecture for $L_p(\chi,\alpha,\psi^{i}\omega, \psi^{j}\omega, 0)$ when $1\leq i<j\leq p-1$.

Let $K_1$ be the compositum of $K$ and $\Q_1$. We computed the minimal polynomial of the Stark unit for $K_1$ over $F$. Let $u_1$ be a root of the minimal polynomial, so $u_1$ is a Stark unit for $K_1$ over $F$. 

In all three examples, the Hecke polynomial of $f$ at $p$ is $x^2 - 1$, so $\alpha = \pm1$ and the geometric Frobenius, $\delta_p$ has order two. For a $\Delta_p$-module $A$ and $a\in A$, let
$$|a|_\alpha = \begin{cases} a\delta_p(a) &\text{ if } \alpha = 1\\
						\frac{a}{\delta_p(a)} &\text{ if } \alpha = -1\end{cases}$$
so $|\cdot|_\alpha:A\rightarrow A^{\delta_p = \alpha}$. (Note that the definition of $|\cdot|_\alpha$ appearing here differs from the one appearing in Section \ref{finishProofKatz} by the scalar $\frac{1}{|\Delta_p|}$.) Let
$$u_{\chi\psi^i,\alpha}^* = \sum_{\sigma\in \Gal(K_1/F)}\chi\psi^i(\sigma)\otimes|\sigma(u_1)|_{\alpha}\in (k_1\otimes\OX_{M_1}^\times)^{(\rho\psi^i)^*, \delta_p = \alpha}.$$

We computed each example to two levels of precision. First to check the results we computed with 60 $p$-adic digits of precision. Then to reproduce and reaffirm the results we computed each example at a higher level of precision. Let $prec$ be the number of $p$-adic digits that each computation was done with. We computed each of the $p$-adic numbers $L_p(\chi,\alpha,\psi^i\omega,\psi^j\omega,0)$ and $\ds\frac{\log_p(u_{\chi\psi^i,\alpha}^*)}{\log_p(u_{\chi\psi^j,\alpha}^*)}$, which lie in the $p$-adic field $\Q_p(\zeta_{p^2})$. The field $\Q_p(\zeta_{p^2})$ has ramification index $p(p-1)$ over $\Q_p$ and was represented in the computer with respect to the uniformizer $\pi = \zeta_{p^2} - 1$. Computing with $prec$ $p$-adic digits in $\Q_p(\zeta_{p^2})$ is $prec\cdot p(p-1)$, $\pi$-adic digits. To verify the conjecture, we calculated the $\pi$-adic valuation of the difference
\begin{equation} \label{numerics}L_p(\chi,\alpha,\psi^i\omega,\psi^j\omega,0) - \frac{\tau(\psi^{-i})}{\tau(\psi^{-j})}\frac{\log_p(u_{\chi\psi^i,\alpha}^*)}{\log_p(u_{\chi\psi^j,\alpha}^*)}.\end{equation}

A number in our computer representation of $\Q_p(\zeta_{p^2})$ is $0$ if it has $\pi$-adic valuation $prec\cdot p(p-1)$. Our data shows that in the cases that we computed, the value of (\ref{numerics}) is extremely close to $0$. The difference in the examples between $(\ref{numerics})$ and $0$ is most likely from rounding error.

\subsection{$F = \Q(\sqrt{17})$, $K = \Q(\sqrt{4 + \sqrt{17}})$, $p=5$}\label{Ex1}

In this example, Conjecture \ref{pAdicStarkAt0} is true because $\rho$ is also the induction of a ray class character $\chi'$ of $F' = \Q(i)$ where $p = 5$ splits (and Conjecture \ref{pAdicStarkAt0} only depends on $\rho$). To see this, define $\chi'$ so that the fixed field of the kernel of $\chi'$ is $K' = \Q(\sqrt{8 + 2i})$. Then since the fixed field of the kernel of $\rho$ is $M = K(\sqrt{4 - \sqrt{17}})$ and we have the relation $\sqrt{4 + \sqrt{17}} + \sqrt{4 - \sqrt{17}} = \sqrt{8 + 2i}$, a simple calculation shows that $\Ind\chi = \rho = \Ind\chi'$. For a further analysis of this situation where there is a ray class character of a real quadratic field and of an imaginary quadratic field where $p$ splits, and such that the induction of the two ray class characters is the same, see chapter 5 of \cite{Ferr}.

We include this example because the units appearing are the Stark units from \ref{StarkUnitDef} associated to the real quadratic field $F = \Q(\sqrt{17})$. 

The character $\psi$ is defined by $\psi(2) = \zeta_5$. Let $a = \ds\frac{1 + \sqrt{17}}{2}$. Then the minimal polynomial of the Stark unit for $K_1/F$ is
\[\niceArrayPolys \begin{array}{c} x^{10} + (-2268731445425a - 3542743970110)x^9 +\\
				 + (101815525268417913200a + 158990319870506526445)x^8 + \\
				+ (-908489137763713280149684575a - 1418653768481195383230297220)x^7 +\\
				+ (1212779745101402982169172133826675a + 1893819622280672026587959027568110)x^6 + \\
				+ (-51814142160111896449580114635979570875a - 80910519433399332983120295909704647352)x^5 +\\
				+ (1212779745101402982169172133826675a + 1893819622280672026587959027568110)x^4 +\\
				+ (-908489137763713280149684575a - 1418653768481195383230297220)x^3 +\\
				+ (101815525268417913200a + 158990319870506526445)x^2 +\\
				+  (-2268731445425a - 3542743970110)x + 1.\end{array}\]
The data for this example is in the following table.

\begin{center}\begin{tabular}{|l|l|l|l|}
\hline
$\alpha$ & (i,j) & $\substack{\pi\text{-adic valuation of (\ref{numerics})}\\ \text{when } prec = 60}$ &$\substack{\pi\text{-adic valuation of (\ref{numerics})}\\ \text{when }prec = 63}$\\
\hline
1 & (1,2) &1141&1260\\
\hline
1 & (1,3) &1140&1260\\
\hline
1 & (1,4) &1140&1261\\
\hline
1 & (2,3) &1140&1260\\
\hline
1 & (2,4) &1140&1260\\
\hline
1 & (3,4) &1141&1260\\
\hline
-1 & (1,2) &1136&1255\\
\hline
-1 & (1,3) &1135&1255\\
\hline
-1 & (1,4) &1135&1255\\
\hline
-1 & (2,3) &1135&1255\\
\hline
-1 & (2,4) &1135&1255\\
\hline
-1 & (3,4) &1135&1257\\
\hline
\end{tabular}\end{center}

\subsection{$F = \Q(\sqrt{-23})$, $K = $ Hilbert class field of $F$, $p = 5$}\label{Ex2}

The character $\psi$ is defined by $\psi(2) = \zeta_5$. The minimal polynomial of the Stark unit for $K_1/F$ is

\[\niceArrayPolys \begin{array}{c}  x^{15} - 832535x^{14} + 65231675x^{13} - 5650639400x^{12} + \\
				+ 15533478425x^{11} - 39376942640x^{10} - 212804236525x^9 - 380541320125x^8 + \\
				- 2607229594750x^7 - 2183192838625x^6 + 3771011381950x^5 - 1207366794625x^4 + \\
				+ 99067277500x^3 - 221569375x^2 + 466875x - 125.\end{array}\]

The data for this example is in the following table.

\begin{center}\begin{tabular}{|l|l|l|l|}
\hline
$\alpha$ & (i,j) & $\substack{\pi\text{-adic valuation of (\ref{numerics})}\\ \text{when }prec = 60}$ &$\substack{\pi\text{-adic valuation of (\ref{numerics})}\\ \text{when }prec = 72}$\\
\hline
-1 & (1,2) &1135&1436\\
\hline
-1 & (1,3) &1135&1436\\
\hline
-1 & (1,4) &1135&1435\\
\hline
-1 & (2,3) &1136&1436\\
\hline
-1 & (2,4) &1135&1435\\
\hline
-1 & (3,4) &1135&1435\\
\hline
\end{tabular}\end{center}

When $\alpha = 1$, we made the same calculation and got for (\ref{numerics}) a $p$-adic number that is not close to $0$. This indicates that when $F$ is imaginary quadratic and $p$ is inert in $F$, the units that appear in Conjecture \ref{pAdicStarkAt0} may not come from the elliptic units from definition \ref{ComplexStarkUnitDef}. For reference we give the first $100$ $\pi$-adic digits of the quantities in (\ref{numerics}) for this example when $\alpha = 1$:

\begin{multline}
\frac{\tau(\psi^{-1})}{\tau(\psi^{-2})}\frac{\log_p(u_{\chi\psi^1,\alpha}^*)}{\log_p(u_{\chi\psi^2,\alpha}^*)} = \\2 + \pi^{5} + 4\pi^{21} + 3\pi^{22} + 3\pi^{23} + 4\pi^{24} + \pi^{25} + 2\pi^{26} + 4\pi^{27} + 4\pi^{28} + 2\pi^{29} + 2\pi^{30} + 3\pi^{31} + \pi^{32} + \pi^{33} +\\ 3\pi^{34} + 3\pi^{35} + \pi^{36} + \pi^{37} + 3\pi^{38} + 3\pi^{39} + 3\pi^{40} + 3\pi^{41} + 4\pi^{42} + 3\pi^{43} + 2\pi^{44} + 2\pi^{45} + 2\pi^{46} + 2\pi^{47} + \pi^{48} + \pi^{49} +\\ \pi^{50} + 2\pi^{52} + 3\pi^{54} + 4\pi^{55} + 3\pi^{56} + \pi^{57} + \pi^{58} + 2\pi^{59} + 2\pi^{61} + 4\pi^{62} + 3\pi^{63} + 2\pi^{64} + 3\pi^{65} + \pi^{66} + 2\pi^{67} + 2\pi^{68} +\\ \pi^{71} + 3\pi^{72} + 2\pi^{73} + 2\pi^{74} + \pi^{75} + 2\pi^{76} + 3\pi^{77} + \pi^{78} + 3\pi^{79} + 2\pi^{80} + \pi^{81} + 2\pi^{82} + 4\pi^{84} + 4\pi^{85} + 2\pi^{86} + 4\pi^{88} +\\ 2\pi^{89} + 3\pi^{90} + 3\pi^{91} + 3\pi^{93} + 2\pi^{94} + 4\pi^{95} + 2\pi^{96} + \pi^{97} + 4\pi^{98} + 2\pi^{100} + O(\pi^{101})
\end{multline}

\begin{multline}
\frac{\tau(\psi^{-1})}{\tau(\psi^{-3})}\frac{\log_p(u_{\chi\psi^1,\alpha}^*)}{\log_p(u_{\chi\psi^3,\alpha}^*)} =\\ 3 + 3\pi^{5} + \pi^{10} + 2\pi^{21} + 4\pi^{22} + 4\pi^{23} + 2\pi^{24} + 3\pi^{25} + 4\pi^{26} + 3\pi^{27} + 3\pi^{28} + 4\pi^{29} + 4\pi^{30} + 3\pi^{31} + \pi^{32} +\\ \pi^{33} + 3\pi^{34} + 4\pi^{36} + \pi^{38} + \pi^{40} + 4\pi^{41} + 4\pi^{42} + 3\pi^{43} + 3\pi^{44} + 4\pi^{45} + 4\pi^{46} + 2\pi^{47} + 2\pi^{48} + 4\pi^{49} + 4\pi^{50} + 3\pi^{51} +\\ 3\pi^{53} + \pi^{54} + 4\pi^{55} + 2\pi^{57} + 2\pi^{58} + 2\pi^{59} + 3\pi^{60} + 3\pi^{62} + 2\pi^{63} + 4\pi^{65} + 3\pi^{66} + 3\pi^{67} + 4\pi^{69} + \pi^{70} + \pi^{71} + 2\pi^{72} +\\ 2\pi^{73} + 2\pi^{74} + 3\pi^{75} + 3\pi^{76} + \pi^{77} + 2\pi^{78} + 4\pi^{80} + 4\pi^{81} + 2\pi^{82} + 3\pi^{84} + 3\pi^{85} + 3\pi^{86} + 2\pi^{87} + 3\pi^{88} + 4\pi^{89} + 4\pi^{90} +\\ 2\pi^{91} + 4\pi^{92} + 4\pi^{93} + \pi^{95} + \pi^{97} + 2\pi^{98} + 2\pi^{99} + 2\pi^{100} + O(\pi^{101})
\end{multline}

\begin{multline}
\frac{\tau(\psi^{-1})}{\tau(\psi^{-4})}\frac{\log_p(u_{\chi\psi^1,\alpha}^*)}{\log_p(u_{\chi\psi^4,\alpha}^*)} =\\ 4 + \pi^{5} + 4\pi^{10} + \pi^{15} + 4\pi^{21} + 3\pi^{22} + 3\pi^{23} + 4\pi^{24} + 4\pi^{26} + 3\pi^{27} + 3\pi^{28} + 4\pi^{29} + \pi^{30} + \pi^{31} + 2\pi^{32} +\\ 2\pi^{33} + \pi^{34} + 3\pi^{35} + \pi^{36} + \pi^{37} + 3\pi^{38} + 3\pi^{39} + 2\pi^{40} + 4\pi^{42} + \pi^{43} + 2\pi^{44} + \pi^{45} + \pi^{46} + 2\pi^{47} + 3\pi^{48} + 2\pi^{49} +\\ \pi^{50} + \pi^{51} + 3\pi^{52} + 4\pi^{53} + 2\pi^{54} + 4\pi^{55} + 3\pi^{57} + 4\pi^{58} + 2\pi^{60} + 4\pi^{62} + 4\pi^{65} + 3\pi^{66} + 3\pi^{67} + \pi^{68} + 4\pi^{69} + 3\pi^{70} +\\ \pi^{71} + \pi^{72} + 4\pi^{73} + 2\pi^{74} + 2\pi^{75} + \pi^{76} + 4\pi^{77} + 4\pi^{80} + \pi^{81} + 3\pi^{82} + \pi^{83} + 2\pi^{85} + 3\pi^{87} + 2\pi^{88} + 2\pi^{89} + 3\pi^{91} +\\ \pi^{92} + \pi^{93} + 2\pi^{95} + \pi^{96} + 2\pi^{97} + 4\pi^{99} + 2\pi^{100} O(\pi^{101})
\end{multline}

\begin{multline}
\frac{\tau(\psi^{-2})}{\tau(\psi^{-3})}\frac{\log_p(u_{\chi\psi^2,\alpha}^*)}{\log_p(u_{\chi\psi^3,\alpha}^*)} =\\ 4 + 2\pi^{5} + 2\pi^{10} + 4\pi^{15} + \pi^{20} + 3\pi^{21} + \pi^{22} + \pi^{23} + 3\pi^{24} + \pi^{25} + 2\pi^{30} + 2\pi^{31} + 4\pi^{32} + 4\pi^{33} + 2\pi^{34} +\\ 3\pi^{35} + 4\pi^{36} + \pi^{37} + 3\pi^{39} + 3\pi^{41} + \pi^{42} + 3\pi^{43} + \pi^{45} + \pi^{46} + 4\pi^{48} + 4\pi^{51} + 3\pi^{53} + 2\pi^{54} + 2\pi^{55} + 3\pi^{56} + \pi^{57} +\\ 2\pi^{58} + 4\pi^{59} + 2\pi^{60} + \pi^{61} + 4\pi^{62} + 3\pi^{63} + 2\pi^{64} + 3\pi^{66} + 3\pi^{67} + \pi^{69} + \pi^{70} + 3\pi^{72} + 2\pi^{73} + \pi^{74} + 4\pi^{75} + 3\pi^{76} +\\ 3\pi^{77} + 2\pi^{78} + 4\pi^{79} + 3\pi^{80} + \pi^{81} + 4\pi^{82} + 3\pi^{83} + 4\pi^{84} + 2\pi^{85} + \pi^{86} + 4\pi^{87} + 4\pi^{88} + 4\pi^{89} + 3\pi^{90} + 2\pi^{91} + 2\pi^{92} +\\ 4\pi^{93} + \pi^{94} + 3\pi^{95} + 3\pi^{96} + 2\pi^{97} + 4\pi^{99} + 2\pi^{100} + O(\pi^{101})
\end{multline}

\begin{multline}
\frac{\tau(\psi^{-2})}{\tau(\psi^{-4})}\frac{\log_p(u_{\chi\psi^2,\alpha}^*)}{\log_p(u_{\chi\psi^4,\alpha}^*)} = \\2 + 2\pi^{5} + \pi^{10} + 3\pi^{21} + \pi^{22} + \pi^{23} + 3\pi^{24} + 2\pi^{25} + 2\pi^{26} + 4\pi^{27} + 4\pi^{28} + 2\pi^{29} + 3\pi^{30} + 3\pi^{35} + 3\pi^{36} +\\ 4\pi^{37} + 3\pi^{38} + 2\pi^{39} + 2\pi^{40} + 2\pi^{41} + 2\pi^{44} + \pi^{45} + 4\pi^{46} + 4\pi^{47} + \pi^{52} + \pi^{53} + \pi^{56} + \pi^{58} + 2\pi^{59} + 2\pi^{64} + 4\pi^{65} +\\ \pi^{68} + 4\pi^{69} + \pi^{70} + 4\pi^{71} + 2\pi^{72} + \pi^{73} + 3\pi^{76} + 4\pi^{77} + 4\pi^{78} + \pi^{79} + 3\pi^{80} + 2\pi^{81} + 2\pi^{82} + 4\pi^{83} + 2\pi^{85} + 2\pi^{86} +\\ 4\pi^{88} + \pi^{89} + \pi^{90} + 3\pi^{91} + 4\pi^{92} + 3\pi^{93} + 4\pi^{94} + 2\pi^{95} + 2\pi^{96} + 2\pi^{97} + 2\pi^{98} + O(\pi^{101})
\end{multline}

\begin{multline}
\frac{\tau(\psi^{-3})}{\tau(\psi^{-4})}\frac{\log_p(u_{\chi\psi^3,\alpha}^*)}{\log_p(u_{\chi\psi^4,\alpha}^*)} = \\3 + 4\pi^{5} + 3\pi^{10} + \pi^{15} + \pi^{21} + 2\pi^{22} + 2\pi^{23} + \pi^{24} + 4\pi^{25} + 2\pi^{26} + 4\pi^{27} + 4\pi^{28} + 2\pi^{29} + \pi^{30} + \pi^{31} +\\ 2\pi^{32} + 2\pi^{33} + \pi^{34} + \pi^{35} + 3\pi^{36} + 2\pi^{37} + \pi^{39} + 3\pi^{40} + 4\pi^{41} + \pi^{42} + 3\pi^{43} + \pi^{44} + 3\pi^{51} + 4\pi^{52} + 4\pi^{54} + \pi^{55} +\\ 4\pi^{56} + \pi^{57} + \pi^{58} + 3\pi^{59} + 4\pi^{60} + \pi^{61} + 4\pi^{62} + 4\pi^{64} + 4\pi^{66} + 2\pi^{67} + 3\pi^{70} + 4\pi^{71} + 2\pi^{72} + 3\pi^{73} + 4\pi^{74} + 4\pi^{75} +\\ 3\pi^{76} + \pi^{77} + \pi^{78} + 3\pi^{79} + 4\pi^{80} + \pi^{82} + 3\pi^{83} + 2\pi^{84} + \pi^{85} + \pi^{87} + \pi^{88} + 4\pi^{89} + 2\pi^{92} + \pi^{93} + 2\pi^{95} + 2\pi^{96} +\\ 2\pi^{97} + 3\pi^{98} + 2\pi^{99} + 4\pi^{100} + O(\pi^{101})
\end{multline}

\begin{multline}
L_p(\chi,\alpha,\psi^1\omega,\psi^2\omega,0) = \\1 + 3\pi^{5} + 3\pi^{10} + \pi^{15} + 2\pi^{21} + 4\pi^{22} + 4\pi^{23} + 2\pi^{24} + 2\pi^{25} + \pi^{30} + 3\pi^{31} + \pi^{32} +\\ \pi^{33} + 3\pi^{34} + 2\pi^{35} + 2\pi^{36} + \pi^{37} + 2\pi^{38} + 3\pi^{39} + 4\pi^{40} + 2\pi^{41} + 3\pi^{42} + \pi^{43} + \pi^{45} + 4\pi^{46} + 2\pi^{47} + 3\pi^{48} + \pi^{50} +\\ 2\pi^{51} + \pi^{52} + 4\pi^{53} + 3\pi^{56} + 3\pi^{58} + 4\pi^{59} + 2\pi^{60} + 3\pi^{61} + 4\pi^{62} + \pi^{63} + 3\pi^{64} + 2\pi^{66} + 3\pi^{68} + 3\pi^{69} + 3\pi^{70} + 4\pi^{73} +\\ 2\pi^{75} + \pi^{76} + 4\pi^{80} + 4\pi^{81} + 3\pi^{82} + 3\pi^{83} + 3\pi^{84} + \pi^{85} + \pi^{86} + \pi^{87} + 2\pi^{88} + 3\pi^{89} + \pi^{91} + 3\pi^{92} + 4\pi^{93} + 4\pi^{94} +\\ \pi^{95} + 2\pi^{96} + 3\pi^{97} + 4\pi^{98} + 3\pi^{99} + \pi^{100} + O(\pi^{101})
\end{multline}

\begin{multline}
L_p(\chi,\alpha,\psi^1\omega,\psi^3\omega,0) = \\1 + \pi^{5} + 4\pi^{21} + 3\pi^{22} + 3\pi^{23} + 4\pi^{24} + 4\pi^{25} + 2\pi^{26} + 4\pi^{27} + 4\pi^{28} + 2\pi^{29} + 4\pi^{30} +\\ 3\pi^{31} + \pi^{32} + \pi^{33} + 3\pi^{34} + 3\pi^{35} + \pi^{36} + \pi^{37} + 3\pi^{38} + 3\pi^{39} + 3\pi^{41} + 4\pi^{42} + 3\pi^{43} + 2\pi^{44} + 3\pi^{47} + 2\pi^{48} + 4\pi^{49} +\\ 3\pi^{50} + 4\pi^{51} + 3\pi^{53} + 2\pi^{54} + 4\pi^{55} + 3\pi^{56} + \pi^{57} + \pi^{58} + 2\pi^{59} + 3\pi^{60} + 2\pi^{61} + 2\pi^{62} + \pi^{64} + 2\pi^{66} + 3\pi^{67} + \pi^{68} +\\ 2\pi^{70} + \pi^{71} + \pi^{72} + 4\pi^{73} + 2\pi^{74} + \pi^{75} + 3\pi^{76} + 4\pi^{78} + \pi^{79} + 4\pi^{82} + 3\pi^{83} + 3\pi^{84} + 4\pi^{85} + 4\pi^{87} + 4\pi^{88} + 3\pi^{89} +\\ 3\pi^{91} + 4\pi^{92} + \pi^{93} + 2\pi^{94} + \pi^{95} + 2\pi^{96} + 4\pi^{97} + 3\pi^{98} + 2\pi^{99} + 2\pi^{100} + O(\pi^{101})
\end{multline}

\begin{multline}
L_p(\chi,\alpha,\psi^1\omega,\psi^4\omega,0) = \\1 + 4\pi^{5} + \pi^{10} + 4\pi^{15} + \pi^{20} + \pi^{21} + 2\pi^{22} + 2\pi^{23} + \pi^{24} + \pi^{25} + \pi^{26} + 2\pi^{27} +\\ 2\pi^{28} + \pi^{29} + 3\pi^{30} + 4\pi^{31} + 3\pi^{32} + 3\pi^{33} + 4\pi^{34} + 2\pi^{35} + \pi^{37} + 4\pi^{38} + 3\pi^{39} + 4\pi^{40} + 3\pi^{43} + 3\pi^{44} + 3\pi^{45} + 4\pi^{46} +\\ 4\pi^{47} + 3\pi^{48} + 3\pi^{49} + 4\pi^{55} + 4\pi^{56} + 4\pi^{57} + 3\pi^{58} + 2\pi^{59} + 2\pi^{60} + 2\pi^{61} + 4\pi^{62} + \pi^{63} + 4\pi^{65} + 2\pi^{66} + \pi^{67} + 2\pi^{68} +\\ 3\pi^{69} + 2\pi^{70} + 2\pi^{71} + 4\pi^{72} + \pi^{73} + 2\pi^{74} + 4\pi^{76} + 3\pi^{77} + \pi^{78} + 4\pi^{79} + 2\pi^{80} + 4\pi^{81} + 2\pi^{82} + 2\pi^{85} + 4\pi^{86} + \pi^{87} +\\ 2\pi^{88} + \pi^{89} + 3\pi^{90} + \pi^{91} + 4\pi^{92} + \pi^{96} + \pi^{97} + \pi^{98} + 4\pi^{99} + 4\pi^{100} + O(\pi^{101})
\end{multline}

\begin{multline}
L_p(\chi,\alpha,\psi^2\omega,\psi^3\omega,0) = \\1 + 3\pi^{5} + 3\pi^{10} + \pi^{15} + 2\pi^{21} + 4\pi^{22} + 4\pi^{23} + 2\pi^{24} + 2\pi^{25} + \pi^{30} + 3\pi^{31} + \pi^{32} +\\ \pi^{33} + 3\pi^{34} + 3\pi^{35} + 2\pi^{36} + \pi^{37} + 2\pi^{38} + 3\pi^{39} + \pi^{40} + 2\pi^{41} + 3\pi^{42} + \pi^{43} + 2\pi^{45} + 4\pi^{46} + 2\pi^{47} + 3\pi^{48} + \pi^{50} +\\ 2\pi^{52} + 3\pi^{54} + \pi^{56} + \pi^{57} + 4\pi^{58} + 2\pi^{59} + 3\pi^{60} + 2\pi^{61} + 2\pi^{62} + 4\pi^{63} + 2\pi^{64} + 4\pi^{65} + 3\pi^{69} + 3\pi^{70} + \pi^{71} + \pi^{73} +\\ 3\pi^{75} + 3\pi^{76} + \pi^{77} + 2\pi^{78} + 4\pi^{79} + 4\pi^{80} + \pi^{82} + \pi^{84} + 4\pi^{87} + 4\pi^{90} + 4\pi^{91} + 3\pi^{92} + 2\pi^{93} + \pi^{94} + \pi^{96} + 2\pi^{97} +\\ 4\pi^{98} + 2\pi^{99} + 3\pi^{100} + O(\pi^{101})
\end{multline}

\begin{multline}
L_p(\chi,\alpha,\psi^2\omega,\psi^4\omega,0) = \\1 + \pi^{5} + 4\pi^{21} + 3\pi^{22} + 3\pi^{23} + 4\pi^{24} + 4\pi^{25} + 2\pi^{26} + 4\pi^{27} + 4\pi^{28} + 2\pi^{29} + 4\pi^{30} +\\ 3\pi^{31} + \pi^{32} + \pi^{33} + 3\pi^{34} + 3\pi^{35} + \pi^{36} + \pi^{37} + 3\pi^{38} + 3\pi^{39} + 3\pi^{41} + 4\pi^{42} + 3\pi^{43} + 2\pi^{44} + 3\pi^{47} + 2\pi^{48} + 4\pi^{49} +\\ 3\pi^{50} + 4\pi^{51} + 3\pi^{53} + 2\pi^{54} + 4\pi^{55} + 3\pi^{56} + \pi^{57} + \pi^{58} + 2\pi^{59} + 3\pi^{60} + 2\pi^{61} + 2\pi^{62} + \pi^{64} + 2\pi^{66} + 3\pi^{67} + \pi^{68} +\\ \pi^{71} + \pi^{72} + 4\pi^{73} + 2\pi^{74} + 3\pi^{75} + 3\pi^{76} + 4\pi^{78} + \pi^{79} + 4\pi^{82} + 3\pi^{83} + 3\pi^{84} + 2\pi^{85} + 3\pi^{86} + \pi^{89} + 4\pi^{90} + 2\pi^{91} +\\ 2\pi^{92} + 4\pi^{93} + \pi^{94} + 2\pi^{95} + 3\pi^{96} + \pi^{97} + 3\pi^{99} + O(\pi^{101})
\end{multline}

\begin{multline}
L_p(\chi,\alpha,\psi^3\omega,\psi^4\omega,0) = \\1 + 3\pi^{5} + 3\pi^{10} + \pi^{15} + 2\pi^{21} + 4\pi^{22} + 4\pi^{23} + 2\pi^{24} + 2\pi^{25} + \pi^{30} + 3\pi^{31} + \pi^{32} +\\ \pi^{33} + 3\pi^{34} + 2\pi^{35} + 2\pi^{36} + \pi^{37} + 2\pi^{38} + 3\pi^{39} + 4\pi^{40} + 2\pi^{41} + 3\pi^{42} + \pi^{43} + \pi^{45} + 4\pi^{46} + 2\pi^{47} + 3\pi^{48} + \pi^{50} +\\ 2\pi^{51} + \pi^{52} + 4\pi^{53} + 3\pi^{56} + 3\pi^{58} + 4\pi^{59} + 2\pi^{60} + 3\pi^{61} + 4\pi^{62} + \pi^{63} + 3\pi^{64} + 2\pi^{66} + 3\pi^{68} + 3\pi^{69} + \pi^{70} + 4\pi^{73} +\\ \pi^{76} + \pi^{80} + 4\pi^{81} + 3\pi^{82} + 3\pi^{83} + 3\pi^{84} + \pi^{85} + 4\pi^{86} + 2\pi^{87} + 3\pi^{88} + \pi^{89} + \pi^{92} + 2\pi^{93} + 3\pi^{94} + 2\pi^{95} + \pi^{96} +\\ \pi^{97} + 2\pi^{98} + 2\pi^{99} + 3\pi^{100} + O(\pi^{101}).
\end{multline}

\subsection{$F = \Q(\sqrt{-31})$, $K = $ Hilbert class filed of $F$, $p = 3$}\label{Ex3}

This example is interesting because it does not satisfy the assumption, $p\nmid [M:\Q]$ (in this example $M = K$). In this example $p = 3$ which divides $[M:\Q] = 6$. The example does satisfy the condition
$$\Delta_1 = \Gal(M_1/\Q)\cong \Gal(M/\Q)\times \Gal(\Q_1/\Q) = \Delta\times\Gamma_1.$$
The character $\psi$ is defined by $\psi(2) = \zeta_3$. The minimal polynomial of the Stark unit for $K_1/F$ is
$$x^9 - 306x^8 - 1143x^7 - 71640x^6 + 60156x^5 + 117180x^4 + 25704x^3 - 7371x^2 + 5022x - 27.$$
The data for this example is in the following table.

\begin{center}\begin{tabular}{|l|l|l|l|}
\hline
$\alpha$ & (i,j) & $\substack{\pi\text{-adic valuation of (\ref{numerics})}\\ \text{when }prec = 60}$ & $\substack{\pi\text{-adic valuation of (\ref{numerics})}\\ \text{when }prec = 77}$\\
\hline
-1 & (1,2) &352&441\\
\hline
\end{tabular}\end{center}

When $\alpha = 1$, as in the previous example, we made the same calculation and got for (\ref{numerics}) a $p$-adic number that is not close to $0$. Again, this indicates that when $F$ is imaginary quadratic and $p$ is inert in $F$, the units that appear in Conjecture \ref{pAdicStarkAt0} may not come from the elliptic units from definition \ref{ComplexStarkUnitDef}. For reference we give the first $100$ $\pi$-adic digits of the quantities in (\ref{numerics}) for this example when $\alpha = 1$:

\begin{multline}
\frac{\tau(\psi^{-1})}{\tau(\psi^{-2})}\frac{\log_p(u_{\chi\psi^1,\alpha}^*)}{\log_p(u_{\chi\psi^2,\alpha}^*)} =\\ 2 + 2*\pi^{3} + \pi^{6} + \pi^{7} + \pi^{8} + \pi^{9} + \pi^{11} + 2*\pi^{13} + 2*\pi^{14} + 2*\pi^{15} + 2*\pi^{16} + \pi^{17} + 2*\pi^{19} + 2*\pi^{20} + \pi^{21} +\\ \pi^{22} + 2*\pi^{24} + \pi^{26} + 2*\pi^{28} + \pi^{29} + 2*\pi^{30} + 2*\pi^{32} + \pi^{33} + \pi^{35} + 2*\pi^{37} + \pi^{40} + 2*\pi^{41} + 2*\pi^{45} + 2*\pi^{46} + 2*\pi^{47} +\\ \pi^{48} + 2*\pi^{49} + \pi^{50} + 2*\pi^{52} + 2*\pi^{54} + \pi^{55} + 2*\pi^{58} + 2*\pi^{59} + \pi^{60} + 2*\pi^{62} + 2*\pi^{63} + \pi^{64} + 2*\pi^{65} + 2*\pi^{66}\\ + \pi^{67} + 2*\pi^{71} + 2*\pi^{74} + 2*\pi^{75} + \pi^{77} + 2*\pi^{78} + \pi^{79} + 2*\pi^{80} + \pi^{81} + 2*\pi^{82} + 2*\pi^{84} + 2*\pi^{87} + \pi^{89} + \pi^{91} + 2*\pi^{92} +\\ 2*\pi^{98} + 2*\pi^{99} + 2*\pi^{100} + O(\pi^{101})
\end{multline}

\begin{multline}
L_p(\chi,\alpha,\psi^1\omega,\psi^2\omega,0) = \\ 1 + 2*\pi^{3} + \pi^{6} + \pi^{7} + \pi^{8} + 2*\pi^{9} + \pi^{11} + 2*\pi^{13} + 2*\pi^{14} + 2*\pi^{16} + \pi^{17} +\\ \pi^{19} + \pi^{20} + 2*\pi^{21} + \pi^{22} + 2*\pi^{24} + \pi^{25} + \pi^{26} + 2*\pi^{27} + 2*\pi^{28} + 2*\pi^{30} +  \pi^{32} + 2*\pi^{33} + 2*\pi^{35} + \pi^{36} + 2*\pi^{38} + \pi^{40} +\\ 2*\pi^{41} + \pi^{45} + 2*\pi^{49} + 2*\pi^{50} + 2*\pi^{51} + 2*\pi^{52} + 2*\pi^{54} + \pi^{55} + \pi^{57} +  \pi^{58} + \pi^{62} + \pi^{63} + \pi^{64} + 2*\pi^{65} + 2*\pi^{68} +\\ 2*\pi^{69} + \pi^{71} + \pi^{74} + \pi^{75} + \pi^{76} + 2*\pi^{78} + 2*\pi^{82} + \pi^{83} + \pi^{84} + 2*\pi^{86} + \pi^{88} + \pi^{90} + 2*\pi^{91} + \pi^{92} + \pi^{94} + \pi^{97} +\\ 2*\pi^{98} + 2*\pi^{99} + \pi^{100} + O(\pi^{101}).
\end{multline}

\bibliography{ThesisToPaperNewEvidence}
\bibliographystyle{abbrv}

\Addresses

\end{document}